\documentclass[12pt,fleqn]{article}
\usepackage{indentfirst,latexsym}
\usepackage{latexsym,amsmath,amssymb}
\usepackage{bm}
\usepackage{graphicx}
\usepackage{amsthm,amsfonts,mathrsfs}
\usepackage{float}
\usepackage{color}
\usepackage[dvipsnames]{xcolor}
\usepackage[colorlinks,citecolor=blue,urlcolor=magenta]{hyperref}
\usepackage{breakurl}
\usepackage{stmaryrd}
\usepackage{bbding}
\usepackage{booktabs}
\usepackage{threeparttable}
\usepackage[figuresright]{rotating}
\usepackage{cite}
\usepackage[top=1.2in, bottom=1.2in, left=1.2in, right=1.2in]{geometry}
\usepackage{algorithm}
\usepackage{algpseudocode}

\numberwithin{equation}{section}
\newtheorem{theorem}{Theorem}[section]

\date{}
\begin{document}
\title{Efficient variants of the CMRH method for solving a sequence of multi-shifted non-Hermitian linear systems simultaneously}
\author{Xian-Ming Gu$^{1,2,}$\footnote{\textit{E-mail address:} \href{guxianming@live.cn}
{guxianming@live.cn};~\href{x.m.gu@rug.nl}{x.m.gu@rug.nl}},\mbox{ }Ting-Zhu
Huang$^{3,}$\footnote{Corresponding author. \textit{E-mail address:}
\href{tingzhuhuang@126.com}{tingzhuhuang@126.com}. Tel.: 86-28-61831016.},\mbox{ }
Bruno Carpentieri$^{4,}$\footnote{Corresponding author. \textit{E-mail address:} \href{bruno.carpentieri@unibz.it}{bruno.carpentieri@unibz.it}. Tel.: +39-047-101-6027}\\ Akira Imakura$^{5,}$\footnote{\textit{E-mail address:}
\href{imakura@cs.tsukuba.ac.jp}
{imakura@cs.tsukuba.ac.jp}},\mbox{ }Ke Zhang$^{6,}$\footnote{\textit{E-mail address:}
\href{xznuzk123@126.com}{xznuzk123@126.com}},\mbox{ }Lei Du$^{7,}$\footnote{\textit{E-mail
address:} \href{dulei@dlut.edu.cn}{dulei@dlut.edu.cn}}\\
{\footnotesize{\it 1. School of Economic Mathematics/Institute of Mathematics,}}\\
{\footnotesize{\it Southwestern University of Finance and Economics, Chengdu, 611130, P.R. China}}\\
{\footnotesize{\it 2. Bernoulli Institute for Mathematics, Computer Science and Artificial Intelligence,}}\\
{\footnotesize{\it University of Groningen, Nijenborgh 9, P.O. Box 407, 9700 AK Groningen, The Netherlands}}\\
{\footnotesize{\it 3. School of Mathematical Sciences,}}\\
{\footnotesize{\it University of Electronic Science and Technology of China, Chengdu,
611731, P.R. China}}\\
{\footnotesize{\it 4. Faculty of Computer Science,~Free University of Bozen-Bolzano}}\\
{\footnotesize{\it Dominikanerplatz 3 - piazza Domenicani, 3 Italy - 39100, Bozen-Bolzano}}\\
{\footnotesize{\it 5. Graduate School of Systems and Information Engineering,}}\\
{\footnotesize{\it University of Tsukuba, 1-1-1 Tennodai, Tsukuba, Ibaraki, 305-8573, Japan}}\\
{\footnotesize{\it 6. College of Arts and Sciences, Shanghai Maritime University,
Shanghai, 201306, P.R. China}}\\
{\footnotesize{\it 7. School of Mathematical Sciences, Dalian University of
Technology, Dalian, 116024, P.R. China}}
}
\maketitle

\begin{abstract}
Multi-shifted linear systems with non-Hermitian coefficient matrices arise in numerical solutions
of time-dependent partial/fractional differential equations (PDEs/FDEs), in control theory, PageRank
problems, and other research fields. We derive efficient variants of the restarted Changing Minimal
Residual method based on the cost-effective Hessenberg procedure (CMRH) for this problem class. Then,
we introduce a flexible variant of the algorithm that allows to use variable preconditioning at each
iteration to further accelerate the convergence of shifted CMRH. We analyse the performance of the new
class of methods in the numerical solution of PDEs and FDEs, also against other multi-shifted Krylov
subspace methods.



\end{abstract}

\section{Introduction}
\label{sec1}
In this paper we introduce efficient iterative methods for the simultaneous solution of a sequence of,
say $t$, shifted non-Hermitian linear systems of the form
\begin{equation}
(A - \sigma_{i} I){\bm x}^{(i)} = {\bm b},\quad\ i = 1,2,\ldots, t,
\label{eq1.1}
\end{equation}
where $A\in\mathbb{C}^{n\times n}$ is a large, sparse and nonsingular matrix, $\sigma_i\in \mathbb{C}$ are $t
$ shifts given at once, $I$ is the $n\times n$ identity matrix, ${\bm x}^{(i)}$ and ${\bm b}$ are solutions
and right-hand side vectors of the $t$ linear systems, respectively. Problem~(\ref{eq1.1}) arises in implicit
numerical solutions of partial differential (PDEs) \cite{VDLKMZ,LNTJAC} and fractional differential equations
(FDEs) \cite{JACWL,RGAFMP}, in control theory \cite{BNDYS,MIADBS}, large-scale eigenvalue computations \cite{TITSUN},
quantum chromodynamics (QCD) applications \cite{JCRTB} and in other computational science problems \cite{AKSTB,GWYCWXQJ,Baumann18}.

Krylov subspace methods are an efficient alternative to sparse direct methods for solving a sequence of multi-shifted linear systems, owing to the shift-invariance property
\begin{eqnarray}
\mathcal{K}_m(A,{\bm b}) = \mathcal{K}_m(A - \sigma_i I,{\bm b}),\quad\ i = 1,2,\ldots,t,
\label{eq1.2}
\end{eqnarray}
where we denote by $\mathcal{K}_m(A,{\bm b}) :=\mathrm{span}\{{\bm b}, A{\bm b},\ldots, A^{m-1}{\bm b}\}$ the Krylov subspace of dimension $m$ generated by $A$ and ${\bm b}$ ~(see~\cite{MIADBS,BJKSS,AFBCG}). By a suitable choice of the initial vectors ${\bm x}^{(i)}_{0}$, for example take ${\bm x}^{(i)}_{0}=0$, the solution of systems~(\ref{eq1.1}) requires a single Krylov basis~\cite{BJKSS,AFBCG}. This approach has shown to effectively reduce storage and algorithmic costs in the analysis of realistic QCD, PageRank and multi-frequency elastic wave propagation problems~\cite{AKSTB,GWYCWXQJ,Baumann18}.

Over the last two decades, several shift-invariant Krylov subspace algorithms have been proposed for solving multi-shifted linear systems with general non-Hermitian coefficient matrices.  Shifted extensions of the restarted generalized minimum residual (GMRES) method~\cite{AFUG,DDRBMW,GGRG,KMSDBS,KMST,Jing2017}, the  restarted full orthogonalization (FOM) method~\cite{VSRF,YFJTZH,JFYGJY} and the restarted Hessenberg method \cite{XMGTZHG} are some relevant examples  built upon the well-known Arnoldi procedure. On the other hand, shifted versions of the quasi-minimal residual (QMR) method and its transpose-free variant (TFQMR)~\cite{RWFS}, the induced dimension reduction (IDR($s$)) ~\cite{LDTSSLZ,MBMBG}
and its QMR form~\cite{MBGGLG}, the biconjugate gradient (BiCG) method and its stabilized and generalized product-type extensions (BiCGStab, BiCGStab($\ell$) and GPBiCG)~\cite{BJKSS,AFBCG,MDRMA}, the biconjugate residual (BiCR) method and its stabilized form (BiCRSTAB)~\cite{XMGTZHJ}  are built upon short-term vector recurrences such as the Bi-Lanczos~\cite[pp. 229-233]{YSaadI} and the $A$-biorthogonalization~\cite[pp. 40-42]{TSEX} procedures. In~\cite{KAESS,KAPBES}, recycling variants of BiCG and BiCGSTAB have been applied to the solution of multi-shifted non-Hermitian linear systems arising in model reduction applications.

In this paper we consider the Changing Minimal Residual method based on the Hessenberg procedure~\cite[pp. 377-382]{JHWK} (shortly, CMRH\footnote{A short note described the
relation between the ELMRES method \cite{GHDS1,DS1999} and Sadok's CMRH method is available at \url{http://ncsu.edu/hpc/Documents/Publications/gary_howell/contents.html\#codes}.}) introduced in \cite{HSCMRH,MHHSA}. The method is based on long-term vector recurrences, like Arnoldi. However, each basis vector ${\bm l}_i$ has $i-1$ components equal to zero and only one component equal to one; each matrix-vector product involved in the computation of ${\bm l}_i$ costs less than the $Nz$ operations required by the Arnoldi procedure,  where $Nz$ is the number of nonzero entries of $A$. Therefore, the method can be cost-effective especially when large Krylov subspaces are generated~\cite{HSCMRH,MHMDHG,GHDS1,DS1999,XMGTZHG}. It has been shown that CMRH and GMRES convergence curves are often comparable~\cite{HSCMRH,MHHSA,HSDBS}. On problems where GMRES exhibits superlinear convergence, so often does CMRH; when GMRES stagnates, CMRH does so as well. We point the reader to \cite{HSDBS,JDTGM} for further discussions, and to~\cite{KZCG,SDMHPM,Lai2011} for some recent developments on the theoretical and algorithmic aspects of the CMRH method~\cite{SDAP,SDMHPM,DS1999}, including efficient parallel implementations~\cite{HSCMRH}.

We propose an extension of (restarted) CMRH for solving multi-shifted linear systems which
preserves the shift-invariance property of Krylov subspaces by forcing the shifted residuals
to be collinear to the~\textit{seed} system residual at every cycle, at moderate extra storage
and arithmetic operations. It is well known that the use of a preconditioner is essential to
accelerate the convergence of Krylov subspace solvers. Many conventional preconditioning
techniques are not suitable for solving shifted linear systems as they do not ensure that
property~(\ref{eq1.2}) holds for the preconditioned systems~\cite{BJKSS,XMGTZHG,MIADBS,MBMBG}.
We present an inner-outer iterative scheme based on nested Krylov subspace methods, where the
inner solver is a multi-shifted Krylov method such as shifted FOM, shifted IDR($s$), shifted
QMRIDR($s$) or shifted BiCGSTAB($\ell$) methods that acts as a preconditioner for an outer
flexible CMRH (FCMRH) solver \cite{KZCG}.

The rest of the current paper is organized as follows. In Section~\ref{sec2}, we briefly review
the (restarted) CMRH method and we extend it to the solution of multi-shifted linear systems.
In Section~\ref{sec3}, we derive a restarted shifted CMRH method that preserves the shift-invariance property of Krylov subspaces by forcing the shifted residuals to be collinear to the~\textit{seed} system residual at every cycle. Implementation details and algorithmic complexity of the new method are discussed. In Section~\ref{sec4}, we propose a cost-effective nested Krylov subspace method based on shifted CMRH for solving multi-shifted linear systems. Section~\ref{sec5} presents numerical evidence of the potential of the new family of methods to solve efficiently shifted linear systems arising from QCD and from the solutions of PDEs/FDEs. In Section~\ref{sec6}, the paper closes with some final remarks.
\section{The CMRH method for shifted linear systems}
\label{sec2}
In this section, we briefly review the restarted CMRH method for solving general non-Hermitian
linear systems $A{\bm x} = {\bm b}$; then, we extend it to the solution of problem~(\ref{eq1.1}).
\subsection{The CMRH method}
The CMRH algorithm for non-Hermitian systems applies the Hessenberg procedure with pivoting\footnote{According to Refs. \cite{GHDS1,DS1999,XMGTZHG}, it cannot be proved that the Hessenberg procedure with pivoting strategy is backward stable in finite precision arithmetic. However, in most of our experiments the backward error is very small.} to compute a basis of the Krylov subspace $\mathcal{K}_m(A,{\bm r}_0)$. Starting from an initial vector ${\bm l}_1 = {\bm r}_0/\alpha$, where $\alpha = ({\bm r}_0)_1$ is the first entry of ${\bm r}_0$, $m$ steps of the Hessenberg procedure yield a basis $\{{\bm l}_1,{\bm l}_2,\ldots,{\bm l}_m\}$ and a matrix decomposition of the form
\begin{equation}
AL_m = L_{m+1}\check{H}_m,
\label{eq2.1}
\end{equation}
where $L_m = ({\bm l}_1,{\bm l}_2,\ldots,{\bm l}_m)$ is an $n\times m$ unit lower trapezoidal matrix and $\check{H}_m\in \mathbb{C}^{(m + 1)\times m}$ is upper Hessenberg.
The CMRH approximation after $m$ steps writes as ${\bm x}_m = {\bm x}_0 + {\bm z}_m$, where ${\bm z}_m \in \mathcal{K}_m(A,{\bm
r}_0)$ solves the following constrained minimization problem
\begin{equation}
\min_{{\bm u}\in\mathbb{C}^{m+1},~{\bm z}\in \mathcal{K}_m(A,{\bm r}_0)}\|{\bm u}\|_2,\quad\mathrm{subject\
to}\ \ A{\bm z} = {\bm r}_0 + L_{m+1}{\bm u}.
\label{eq2.2}
\end{equation}
Since the columns of $L_{m+1}$ are not orthogonal, vector ${\bm z}_m$ cannot be computed by minimizing directly the residual norm $\|{\bm u}\|_2$  in Eq.~(\ref{eq2.2}). Instead, CMRH computes ${\bm z}_m$ by minimizing the quasi-residual norm
\begin{equation}
\min_{{\bm z}\in \mathcal{K}_m(A,{\bm r}_0)}\|L^{\dag}_{m+1}(A{\bm z}-{\bm r}_0)\|_2,
\label{eqx.q}
\end{equation}
similarly to the well-known QMR method \cite[pp. 236-241]{YSaadI}. In Eq.~(\ref{eqx.q}), we denote by $L^{\dag}_{m+1}$  the pseudo-inverse of $L_{m+1}$; however, note that any left inverse of $L_{m+1}$ would work here~\cite{HSDBS}. Problem~(\ref{eqx.q}) can be interpreted as a standard residual minimization using a semi-norm~\cite{HSCMRH,MHHSA}. By writing ${\bm z} = L_m{\bm y} \in \mathcal{K}_m(A,{\bm r}_0)$, Eq.~(\ref{eqx.q}) is implemented in CMRH as the following small least squares problem
\begin{equation}
\begin{split}
\min_{{\bm y}\in\mathbb{C}^m}\|L^{\dag}_{m+1}({\bm r}_0 - AL_m{\bm y})\|_2& = \min_{{\bm y}\in\mathbb{C}^m}\|L^{\dag
}_{m+1}(\alpha L_{m+1}{\bm e}^{(k+1)}_1 - L_{m+1}\check{H}_m {\bm y})\|_2    \\
& = \min_{{\bm y}\in\mathbb{C}^m}\|\alpha {\bm e}^{(k+1)}_1 - \check{H}_m {\bm y}\|_2.
\end{split}
\label{eq2.3}
\end{equation}

A detailed theoretical analysis of the method can be found in~ \cite{HSCMRH,MHHSA}. In particular, the result below compares  the residual norms computed after $m$ iterations of the CMRH and the GMRES methods:
\begin{theorem} (\cite[Theorem 4]{HSCMRH} and \cite[Theorem 1]{HSDBS})
Let ${\bm r}^{G}_m$ and ${\bm r}^{c}_m$ be the GMRES and CMRH residuals at the $m$th iteration
beginning with the same initial residual ${\bm r}_0$ (e.g., ${\bm r}_0 = {\bm b}$), respectively.
Then
\begin{equation*}
||{\bm r}^{G}_m||_2 \leq ||{\bm r}^{c}_m||_2 \leq \kappa(L_{m +1})||{\bm r}^{G}_m||_2,
\end{equation*}
where $\kappa(L_{m + 1}) = ||L_{m+1}||_2||L^{\dag}_{m+1}||_2$ is the condition number of $L_{m + 1}$.
\label{thm1}
\end{theorem}

Duintjer Tebbens and Meurant have proved that any non-increasing residual norm history is possible for the CMRH method with any set of nonzero eigenvalues of the system matrix~\cite{JDTGM}. Therefore, the distribution of the eigenvalues alone may not play any role in the convergence. However, similarly to other Krylov methods, for many problems and applications a tightly clustered spectrum around a single point away from the origin is favourable to achieve fast convergence, whereas widely spread eigenvalues and/or clusters close to zero are often bad. Being based on long-term vector recurrences, like GMRES also CMRH may need to be restarted to control the growing costs of the Hessenberg procedure. To date, much fewer results are available on the convergence of the restarted CMRH algorithm (referred to as CMRH($m$)).
\subsection{The CMRH method for shifted systems}
\label{sec2.2}
The starting point to develop a shift-invariant extension of the CMRH method for solving
a sequence of multi-shifted linear systems is the shifted Hessenberg relation
\begin{equation}
(A - \sigma_iI)L_m = L_{m + 1}\check{H}_m(\sigma_i),\quad\ \check{H}_m(\sigma_i):= \check{H}_m(\sigma_{seed})
- \sigma_i
\begin{bmatrix}
{I_m}\\
{\bm 0}\\
\end{bmatrix}
\in \mathbb{C}^{(m + 1)\times m},
\label{eq3.1}
\end{equation}
where $I_m$ is the identity matrix of order $m$, and $L_m, \check{H}_m(\sigma_{seed})$ are factors of
the matrix decomposition~(\ref{eq2.1}) independent of $\sigma_i$. By no lack of generality we can assume
that the shift of \textit{seed} system is zero, i.e. $\sigma_{seed} = 0$; otherwise, if $\sigma_{seed}
\neq 0$, we can rewrite Eq.~(\ref{eq1.1}) for $A:= A - \sigma_{seed}I$ and $
\sigma_i :=\sigma_i - \sigma_{seed}$.

By using relation (\ref{eq3.1}), the following shift-dependent CMRH quasi-minimization problem equivalent to Eq.~(\ref{eq2.3}) is derived
\begin{equation}
\min_{{\bm y}\in \mathbb{C}^m}\|L^{\dag}_{m + 1}{\bm r}^{(i)}_m\|_2 = \min_{{\bm y}\in \mathbb{C}^m}\|\alpha
{\bm e}_1 - \check{H}_m(\sigma_i){\bm y}\|_2.
\label{eq3.2x}
\end{equation}

The small-size least squares problem (\ref{eq3.2x}) can be solved in $\mathcal{O}(m^2)$ operations for each
shift $\sigma_i$. Therefore, the most time-consuming part of the algorithm remains the construction of the
Krylov basis $L_m$, but this has to be performed only once. Since the initial residuals must be shift independent,
${\bm x}_0 = {\bm 0}$ should be used as initial vector for all shifted systems. Then, the whole sequence~(\ref{eq1.1})
can be solved simultaneously without additional matrix-vector products in terms of the \textit{seed} system.
\section{The restarted shifted CMRH method}
\label{sec3}
In this section, we present the restarted version of the shifted CMRH algorithm, and a convergence analysis.
We conclude the section with a complexity study of the new method compared to the restarted shifted GMRES solver.
\subsection{The restarted shifted CMRH method}
\label{sec3.1}
Shifted CMRH developed in Section~\ref{sec2.2} may suffer from memory problems in the case long vector recurrences
are generated by the Hessenberg procedure, similarly to its unshifted counterpart. To alleviate such costs, it may
be necessary to restart the algorithm after every, say $m$, Hessenberg steps.  Upon restarting, the residual vectors
${\bm r}^{(i)}_m$ obtained from the quasi-minimum residual condition are not collinear in general, and therefore the
shift-invariance property~(\ref{eq1.2}) may not be maintained. We impose on the \textit{add} systems the collinearity
condition used by Frommer and Gl\"{a}ssner in the restarted shifted GMRES method~\cite{AFUG}, namely
\begin{equation}
{\bm r}^{(i)}_m = \gamma^{(i)}_m {\bm r}_m,\quad\ \gamma^{(i)}_m\in\mathbb{C},
\label{eq3.2}
\end{equation}
where ${\bm r}_m := {\bm r}^{(\sigma_{seed})}_m$ is the residual vector of \textit{seed} system, to ensure that
property~(\ref{eq1.2}) continues to hold at restart. For \textit{seed} system, however, the same quasi-minimum residual
condition enforced on the residual vector by the restarted shifted CMRH method is used. The following
minimization problem is solved for the seed system
\begin{equation*}
{\bm y}_m = \arg \min_{{\bm y}\in \mathbb{R}^m}\|\alpha {\bm e}_1 - \check{H}_m{\bm y}\|_2,
\end{equation*}
where $\check{H}_m := \check{H}_m(\sigma_{seed})$ from Eq.~(\ref{eq3.1}) and ${\bm y}_m:= {\bm y}^{(\sigma_{seed})}_m$.
The residual vector ${\bm r}_m$ appearing in Eq.~(\ref{eq3.2}) can be written as
\begin{equation*}
{\bm r}_m = L_{m + 1}{\bm u}_{m + 1},\quad\ {\bm u}_{m + 1}:= \alpha{\bm e}_1 - \check{H}_m{\bm y}_m\in\mathbb{R}^{m + 1}.
\end{equation*}
Then, from the collinearity condition (\ref{eq3.2}), it follows for the \textit{add} systems
\begin{equation*}
\begin{array}{lllll}
{\bm r}^{(i)}_m =\gamma^{(i)}_m {\bm r}_m\qquad && \Leftrightarrow &&\qquad  {\bm b} - A(\sigma_i)({\bm x
}^{(i)}_0 + L_m{\bm y}^{(i)}_m) = \gamma^{(i)}_m L_{m+1}{\bm u}_{m+1}\\
\qquad && \Leftrightarrow &&\qquad  {\bm r}^{(i)}_0 - A(\sigma_i) L_m{\bm y}^{(i)}_m = L_{m+1}{\bm u}_{m+1}
\gamma^{(i)}_m\\
\qquad && \Leftrightarrow &&\qquad  \gamma^{(i)}_0 {\bm r}_0 - L_{m+1}\check{H}_m(\sigma_i){\bm y}^{(i)}_m = L_
{m+1}{\bm u}_{m+1}\gamma^{(i)}_m\\
\qquad && \Leftrightarrow &&\qquad  L_{m+1}(\check{H}_m(\sigma_i){\bm y}^{(i)}_m + {\bm u}_{m+1}\gamma^{(i)
}_m) = \gamma^{(i)}_0{\bm r}_0\\
\qquad && \Leftrightarrow &&\qquad  \check{H}_m(\sigma_i){\bm y}^{(i)}_m + {\bm u}_{m+1}\gamma^{(i)
}_m = \gamma^{(i)}_0\alpha{\bm e}_1\\
\end{array}
\end{equation*}
using the matrix representation of the Hessenberg procedure; see \cite{AFUG} for details.

From the above relation, ${\bm y}^{(i)}_m$ and $\gamma^{(i)}_m$ can be read as solutions of the $(m+1) \times (m+1)$ linear systems
\begin{equation}
\begin{bmatrix}
\check{H}_m(\sigma_i)  &\mid  & {\bm u}_{m + 1}
\end{bmatrix}
\begin{bmatrix}
{\bm y}^{(i)}_m \\
\gamma^{(i)}_m
\end{bmatrix} = \gamma^{(i)}_0\alpha {\bm e}_1,
\label{eq3.3}
\end{equation}
for $i = 1,2,\ldots,t_s$. Systems (\ref{eq3.3}) are upper Hessenberg, and are solved efficiently using Givens rotations.
The complete restarted shifted CMRH method is summarized in Algorithm~\ref{alg2}.
\begin{algorithm}[!htpb]
\caption{The restarted shifted CMRH method with pivoting.}
\begin{algorithmic}[1]
\Require
the coefficient matrix $A$ (or a user-defined function that applies $A$ to a vector); the right-hand side vector ${\bm b}$;
the set of shifts $\{\sigma_i\}_{i=1,\ldots,t}$; the dimension $m$ of the Krylov subspace; the maximum number of
outer iterations, $maxit$.
\Ensure
the set of solution vectors ${\bm x}^{(i)}$ of the sequence of multi-shifted linear systems.
  \State Choose the restart dimension $m$ and the initial guess ${\bm x}_0, {\bm x}^{(i)}_0$ such that
         ${\bm r}^{(i)}_0 = \gamma^{(i)}_0 {\bm r}_0$, e.g. ${\bm x}_0 = {\bm x}^{(i)}_0 = {\bm 0}$ (also
         implies that $\gamma^{(i)}_0 = 1$)
  \State Set ${\bm q} = [1,2,\ldots n]^T$ and determine $j_0$ such that $|({\bm r}_0)_{j_0}|
         = \|{\bm r}_0\|_{\infty}$
  \State Set $\alpha = ({\bm r}_0)_{j_0},\ {\bm l}_1 = {\bm r}_0/\alpha$ and $({\bm q})_1 \leftrightarrow
         ({\bm q})_{j_0}$, where $\leftrightarrow$ is used to swap contents.
  \For{$j = 1,2,\ldots,m$,}
  \State Compute ${\bm u} = A{\bm l}_j$
  \For{$k = 1,2,\ldots,j$,}
  \State $\check{h}_{k,j} = ({\bm u})_{(\bm q)_k}$
  \State ${\bm u} = {\bm u} - \check{h}_{k,j}{\bm l}_k$
  \EndFor
  \If{$j < n$ and ${\bm u}\neq {\bm 0}$}
  \State Determine $j_{0}\in \{j+1,\ldots,n\}$ such that $ |({\bm u})_{({\bm q})_{j_0}}| = \|({\bm
         u})_{({\bm q})_{j + 1}:({\bm q})_n}\|_{\infty}$
  \State $\check{h}_{j + 1,j} = ({\bm u})_{({\bm q})_{j_0}}$,\ \ ${\bm l}_{j+1} = {\bm u}/\check{
         h}_{j + 1,j};\ ({\bm q})_{j + 1} \leftrightarrow ({\bm q})_{j_0}$
  \Else
  \State $\check{h}_{j + 1,j} = 0$; Stop
  \EndIf
  \EndFor
  \State Define the $(m + 1)\times m$ Hessenberg matrix $\hat{H}_m = \{\check{h}_{k,j}\}_{1\leq k\leq m+1, 1\leq
         j \leq m}$ and collect the matrix $L_m = [{\bm l}_1,{\bm l}_2,\ldots,{\bm l}_m]$
  \State Compute ${\bm y}_m = \arg\min\limits_{{\bm y}\in \mathbb{C}^m} = \|\alpha{\bm e}_1 - \hat{H}_m{\bm y}\|_2$
         and set ${\bm u}_{m + 1} = \alpha {\bm e}_1 - \hat{H}_m {\bm y}_m$
  \State ${\bm x}_m = {\bm x}_0 + L_m{\bm y}_m,\ {\bm r}_m = {\bm r}_0 - AL_m{\bm y}_m = L_{m+1}{\bm u}_{m+1}$
  \For{$i = 1,2,\ldots,t_s$}
  \State Solve $
\begin{bmatrix}
\check{H}_m(\sigma_i) \mid {\bm u}_{m+1}
\end{bmatrix}
\begin{bmatrix}
{\bm y}^{(i)}_m \\
\gamma^{(i)}_m
\end{bmatrix} = \gamma^{(i)}_0\alpha{\bm e}_1$
\State ${\bm x}^{(i)}_m = {\bm x}^{(i)}_0 + L_m {\bm y}^{(i)}_m,~{\bm r}^{(i)}_m = {\bm r}^{(i)}_0 - AL_m {\bm
       y}^{(i)}_m = \gamma^{(i)}_m{\bm r}_m$
  \EndFor
  \State Restart: if converged then stop; otherwise update ${\bm x}_0 := {\bm x}_m,{\bm x}^{(i)}_0 := {\bm x}^{(i)}_m,
  {\bm r}_0 := {\bm r}_m,{\bm r}^{(i)}_0 := {\bm r}^{(i)}_m,\gamma^{(i)}_0: = \gamma^{(i)}_m$ and goto step 2-3.
\end{algorithmic}
\label{alg2}
\end{algorithm}

The following convergence result can be given for the restarted shifted CMRH method.
\begin{theorem}
If the coefficient matrix $A$ is positive definite, i.e., $(A{\bm x},{\bm x}) > 0, \forall{\bm x} \neq {\bm
0}$, and $\sigma_i < 0$ for all $i$, and if the restarted CMRH method converges for the seed system,
then the restarted shifted CMRH method also converges for the \textit{add} systems for every restart frequency $m$.
Moreover, we have
\begin{equation}
\|{\bm r}^{(i)}_{j,m}\|_2 \leq |\gamma^{(i)}_0|\cdot\|{\bm r}_{j,m}\|_2
\label{eq3.4}
\end{equation}
for all $j$, where ${\bm r}^{(i)}_{j,m}$ and ${\bm r}_{j,m}$ represent the $j$-th residual vectors for
the add and the seed systems in a generic cycle, respectively.
\label{them3.1}
\end{theorem}
\noindent\textbf{Proof}. Let $\zeta$ be any zero of the CMRH residual polynomial~\cite{Lai2011} after $m$ iterations, $\varphi_m$. Since $\varphi_m(0) = 1$, it is $\zeta \neq
0$ and we can write
\begin{equation*}
\varphi_m(t) = \Big(1 - \frac{t}{\zeta}\Big)\tau_{m - 1}(t),
\end{equation*}
with $\tau_{m-1}$ a polynomial of degree $m-1$ such that $\tau_{m-1}(0) = 1$.
Denote ${\bm w} = \tau_{m-1}(A){\bm r}_0$ and ${\bm r} = A{\bm w}$. Then
\begin{equation*}
\|\varphi_m(A){\bm r}_0\|_2 = \Big\|{\bm w} - \frac{1}{\zeta}{\bm r}\Big\|_2.
\end{equation*}
For $\gamma\in \mathbb{C}$, the functional $\|{\bm w} - \gamma {\bm r}\|_2$ is minimized for
\begin{equation*}
\gamma^{\star} = \frac{({\bm r},{\bm w})}{({\bm r},{\bm r})}.
\end{equation*}
Assuming that $\|\varphi_m(A){\bm r}_0\|_2$ has been minimized, the norm $\|L^{\dag}_{m+1}
{\bm r}_m\|_2$ is also minimal due to the inequality $\|L^{\dag}_{m+1}{\bm r}_m\|_2 \leq \|L^{\dag
}_{m+1}\|_2\|\varphi_m(A){\bm r}_0\|_2$. We conclude that
\begin{equation*}
\frac{1}{\zeta} = \gamma^{\star} = \frac{({\bm r},{\bm w})}{({\bm r},{\bm r})} = \frac{({\bm r},
A^{-1}{\bm r})}{({\bm r},{\bm r})}
\end{equation*}
and $({\bm r}, A^{-1}{\bm r})/({\bm r},{\bm r}) = (A^{-H}{\bm r}/\|{\bm r}\|_2,{\bm r}/\|{\bm r}\|_2)
\in F(A^{-H})$. Here we denote by $F(A)$ the field of values
$F(A)= \left\{(A{\bm x}, {\bm x})|{\bm x}\in\mathbb{C}^n, \|{\bm x}\|_2 = 1\right\}\subset \mathbb{C}.$

On the other hand, since $A$ is positive definite and $(A{\bm x}, {\bm x}) =  (A^{-H}{\bm y},{\bm y})$ for ${\bm y} = A{\bm x}$, $F(A^{-H})$ is also contained in the right half-plane like $F(A)$. In conclusion, it holds ${\rm Re}\left(\frac{1}{\zeta}\right) > 0$ and thus ${\rm Re}(\zeta) > 0$. Since $\sigma_i < 0$, it is $\varphi(\sigma_i) \neq 0$. Similarly to~\cite[Lemma 2.4]{AFUG} for the restarted shifted GMRES method, this condition ensures that Eq.~(\ref{eq3.3}) has a unique solution. It follows that the restarted shifted CMRH method converges for both \textit{seed} and \textit{add} systems for each restart frequency $m$. Finally, Eq.~(\ref{eq3.4}) can be proved following a similar argument to \cite[Theorem 3.3]{AFUG}. \hfill $\Box$

\textbf{Remark~1}. Analogously to the restarted shifted GMRES method, since $ \gamma^{(i)}_0 = 1$ from ${\bm x}_0 = {\bm
x}^{(i)}_0 = {\bm 0}$, according to Eq.~(\ref{eq3.4}) the \textit{add} systems may converge more rapidly than
the \textit{seed} system. If not, the shift switching technique\footnote{In general, the \textit{seed} system may converge
faster than the \textit{add} systems. In this case, for efficient computation, it needs to change the \textit{seed} system
into one of the rest systems, this strategy is called the seed switching technique.}\cite{TSTHSLZ} can be applied.

\textbf{Remark~2}. In Theorem \ref{them3.1}, restarted CMRH is supposed to be convergent on the positive
definite \textit{seed} system. Unfortunately, there exist very few results in the literature on the
convergence of restarted CMRH, even in the positive definite case; the topic is largely unexplored~\cite{Meurant17,HSCMRH,DS1999,HSDBS,JDTGM}.
According to our computational experience, restarted shifted CMRH often enjoys similar convergence
behaviour to the more costly restarted shifted GMRES method; in our experiments (see Section~\ref{sec5}),
in some cases it could even handle certain shifted systems where the latter failed.
\subsection{Computational cost of the restarted shifted CMRH method}
\label{sec3.2}
Following the idea introduced in our recent work \cite{XMGTZHG,Jing2017}, in Table~\ref{tabku1} we present
a comparative  complexity analysis between shifted CMRH($m$), shifted GMRES($m$) and shifted SGMRES($m$)
\cite{Jing2017} at equal restarting frequency $m$. The main difference in terms of floating-point
operations (FLOPs) between shifted SGMRES($m$) and shifted GMRES($m$) is due to the cost of applying the Givens
rotations in the least-squares solution for \textit{seed} system. Shifted SGMRES($m$) often requires less FLOPs than
shifted GMRES($m$)~\cite{Jing2017}. On the other hand, shifted CMRH($m$) may require less FLOPS than both shifted SGMRES($m$) and shifted GMRES($m$) that are based on the more costly Arnoldi procedure.
\begin{sidewaystable}[!htpb]\tabcolsep=4pt
\begin{center}
\begin{threeparttable}
\caption{Computational cost of a generic cycle of the shifted versions of GMRES, SGMRES and CMRH}
\begin{tabular}{clll}
\hline
Index   & Shifted GMRES & Shifted SGMRES          & Shifted CMRH \\
\hline
$\aleph_1$  & $2mNz + (2m^2 + 5m + 3)n$  & $2mNz + (2m^2 + 5m + 3)n$    & $2mNz + (m + 2)n + m(m+1)(n - \frac{m-1}{3})$ \\
$\aleph_2$  & FLOPs with an Hessenberg matrix $H^{G}_m$ & FLOPs with an upper triangular matrix $R_m$ & FLOPs with an Hessenberg matrix $H^{c}_m$\\
$\aleph_3$  & $t_s(2mn)$        & $t_s(2mn)$           & $t_s(2mn)$ \\
$\aleph_4$  & FLOPs with $\begin{bmatrix}\check{H}^{G}_m(\sigma_i) \mid {\bm u}_{m+1}\end{bmatrix}$
& FLOPs with $R_m - \sigma_i[V^{H}_m{\bm r}_{0,1}, I_{m-1}]$ & FLOPs with $\begin{bmatrix}
\check{H}^{c}_m(\sigma_i) \mid {\bm u}_{m+1}\end{bmatrix}$\\
\hline
\end{tabular}
\label{tabku1}
\begin{tablenotes}
\item $\aleph_1$ represents the cost of iterates for \textit{seed} system;
\item $\aleph_2$ means the cost of least-squares solve for \textit{seed} system;
\item $\aleph_3$ represents the cost of necessary vector updates;
\item $\aleph_4$ means the cost of least-squares solve for $t_s - 1$ add systems ($i = 1,\ldots,t_s - 1$).
\item \textbf{Remark}: In fact, it is useful to mention that the $\mathcal{O}(m^2)$ operations
are needed for coping with the problems in $\aleph_2$ and $\aleph_4$ of those above three solvers, respectively.
However, to solve an upper-triangular least-squares problem (or linear system) is still slightly cheaper
than to solve an upper-Hessenberg least-squares problem (or linear system).
\end{tablenotes}
\end{threeparttable}
\end{center}
\end{sidewaystable}
\section{Inner-outer variants of CMRH for shifted systems}
\label{sec4}
It is celebrated that preconditioning is an essential ingredient to accelerate the convergence of Krylov subspace methods, including their shifted variants.
Without further assumptions on the preconditioners $M(\sigma_i)$ applied to the $t$ linear systems~(\ref{eq1.1}), the shift-invariant property~(\ref{eq1.2}) may not be preserved for the preconditioned Krylov subspaces~\cite{MIADBS,BJKSS,MBMBG}. One would like to find a matrix $M$, independent of $\sigma_i$, satisfying
\begin{equation}
\mathcal{K}_m(AM^{-1}, {\bm b}) = \mathcal{K}_m((A - \sigma_i I)M(\sigma_i)^{-1}, {\bm b}),
\label{eq4.1}
\end{equation}
so that property~(\ref{eq1.2}) would hold true in the preconditioned case as well. It is not hard to find that Eq.~(\ref{eq4.1}) is satisfied if the preconditioned shifted matrix can be written as a shifted preconditioned matrix, that is
\begin{equation}
(A - \sigma_i I)M(\sigma_i)^{-1} =AM^{-1} - \eta_{i}I.
\label{eq4.2}
\end{equation}
Although $M$ is not needed to compute a basis of $\mathcal{K}_m((A - \sigma_i I)M(\sigma_i)^{-1}, {\bm b})$, it must be known explicitly to compute the solutions ${\bm x}^{(i)}$ of the unpreconditioned system from the knowledge of the solution ${\bm y}^{(i)}$ of the right-preconditioned system; see e.g., \cite{GGXZLL} for details.

Inspired by~\cite{AKSTB,GGXZLL}, in this work we use flexible preconditioning to solve sequence~(\ref{eq1.1}). Flexible preconditioning means that a different preconditioner can be applied at every iteration $j$ of an iterative Krylov method, see e.g. \cite{YSaadA,VSDBSF,BCIDLGGS} and our recent work \cite{KZCG,DSTHBCYJ}. If different preconditioners $P_j,P_j(\sigma)$ are used at every iteration $j$, a relation similar to Eq.~(\ref{eq4.2}), namely
\begin{equation}
(A - \sigma_i I)P_j(\sigma_i)^{-1} = \alpha_j(\sigma_i)AP^{-1}_j - \beta_j(\sigma_i)I,
\label{eq1.23}
\end{equation}
must hold to ensure the shift invariance property of the preconditioned Krylov subspace given by Eq.~(\ref{eq4.1}). In Eq.~(\ref{eq1.23}), $\alpha_j$ and $\beta_j$ are parameters dependent on the shifts $\sigma_i$. At each iteration $j$\emph{}, the preconditioner is applied to a vector ${\bm v}_j$ in the form
\begin{equation}
(A - \sigma_i I)P_j(\sigma_i)^{-1}{\bm v}_j = \alpha_j(\sigma_i)AP^{-1}_j {\bm v}_j - \beta_j(\sigma_i){\bm v}_j.
\label{eq1.24}
\end{equation}
Next, we determine conditions on the coefficients $\alpha_j$'s and $\beta_j$'s to ensure that Eq.~(\ref{eq1.23}) and Eq.~(\ref{eq1.24}) hold. Preconditioning is applied by using an inner solver in a multi-shift Krylov  method, and the preconditioned vectors
\begin{equation*}
{\bm z}_j = P^{-1}_j{\bm v}_j,\qquad\ \ {\bm z}^{(i)}_j = P_j(\sigma_i)^{-1}{\bm v}_j,
\end{equation*}
are computed via a truncated multi-shifted Krylov subspace solver. Therefore, the corresponding (inner) residuals are given by
\begin{eqnarray}
{\bm r}_j& = &  {\bm v}_j - A{\bm z}_j = {\bm v}_j - AP^{-1}_j {\bm v}_j, \label{eq1.25}\\
{\bm r}^{(i)}_j &= & {\bm v}_j - (A - \sigma_i I){\bm z}^{(i)}_j = {\bm v}_j -
(A - \sigma_i I)P_j(\sigma_i)^{-1}{\bm v}_j.  \label{eq1.26}
\end{eqnarray}
We require the residuals (\ref{eq1.25})-(\ref{eq1.26}) of the inner method to be collinear, i.e.
\begin{equation}
\exists~\gamma^{(i)}_j \in\mathbb{C}:\quad\ \gamma^{(i)}_j{\bm r}_j = {\bm r}^{(i)}_j.
\label{eq1.27}
\end{equation}
Note that the collinearity factors $\gamma^{(i)}_j$ change at each iteration $j$, for every shift $\sigma_i$.  Eq.~(\ref{eq1.27}) is satisfied automatically by methods such as shifted FOM \cite{VSRF,YFJTZH,JFYGJY}, shifted BiCGStab \cite{BJKSS,AFBCG}, shifted BiCRStab \cite{XMGTZHJ}, shifted IDR($s$) \cite{MBMBG,LDTSSLZ} and by the restarted shifted Hessenberg method~\cite{XMGTZHG}. From the conclusion in \cite{MBMBG},
$\alpha_j$'s and $\beta_j$'s can be determined using~Eq.~(\ref{eq1.24}) and the collinearity relation~(\ref{eq1.27}). The residuals are collinear if
\begin{equation*}
a_j = \gamma^{(i)}_j,\qquad \beta_j = \alpha_j - 1 = \gamma^{(i)}_j - 1
\end{equation*}
at every (outer) iteration $1\leq j \leq m$. One can show that the following relation,
\begin{equation}
(A - \sigma_iI){\bm z}^{(i)}_j = \gamma^{(i)}_j A{\bm z}_j -(\gamma^{(i)}_j - 1){\bm v}_j,
\label{eq1.28}
\end{equation}
holds or, in terms of the flexible preconditioners $P_j$ and $P_j(\sigma_i)$,
\begin{equation*}
(A - \sigma_iI)P_j(\sigma_i)^{-1}{\bm v}_j = \Big(\gamma^{(i)}_j AP^{-1}_j - (\gamma^{(i)}_j
-1)I\Big){\bm v}_j,\quad\  1\leq j \leq m.
\end{equation*}
It remarks that $\alpha_j$'s and $\beta_j$'s do depend on $\sigma_i$, since the collinearity factors $\gamma^{(i)}_j$ change for each shift.

Based on the strategy proposed by Baumann and van Gijzen in \cite{MBMBG}, we now present a new nested multi-shifted Krylov solver in which the shifted Hessenberg (\texttt{msHessen}) method~\cite{XMGTZHG} is used as an inner preconditioner and flexible CMRH (\texttt{FCMRH}) is used as the outer Krylov iteration. We note that shifted FOM and flexible shifted GMRES  \cite{YSaadA} are related to the shifted Hessenberg and the flexible shifted CMRH methods, respectively.
The Hessenberg relation given by Eq.~(\ref{eq2.1}) can be extended as follows,
\begin{equation*}
AZ_m = L_{m+1}\check{H}_m,\quad\quad (A - \sigma_i I)Z^{(\sigma_i)}_{m} = L_{m+1}
\underline{H}_m(\sigma_i),
\end{equation*}
where at iteration $1\leq j \leq m$ flexible preconditioning is applied in the form $z_j = P^{-1}_j{\bm v}_j$,
$z^{(\sigma_i)}_j = P^{-1}_j(\sigma_i){\bm v}_j$, with $Z_m = [{\bm z}_1,\ldots, {\bm z}_m]$ and $Z^{(\sigma_i)}_m = [{\bm z}^{(\sigma_i)}_1,\ldots, {\bm z}^{(\sigma_i)}_m]$. It follows that both the Hessenberg process and the strategy proposed in \cite{MBMBG} yield the modified Hessenberg
matrix
\begin{equation*}
\underline{H}_m(\sigma_i) = (\check{H}_m - \underline{I}_m)\Gamma^{(i)}_m + \underline{I}_m,
\end{equation*}
where $\underline{I}_m$ is the $m \times m$ identity matrix with an extra column of zeros appended.
The consecutive collinearity factors of the inner method then appear on a diagonal matrix $\Gamma^{(i)}_m$ defined as
\begin{equation}
\Gamma^{(i)}_m \equiv
\begin{pmatrix}
\gamma^{(i)}_1 \\
&\gamma^{(i)}_2\\
&&\ddots\\
&&&\gamma^{(i)}_m\\
\end{pmatrix}\in \mathbb{C}^{m\times m}.
\label{eq4.4}
\end{equation}

After $m$ outer iterations of the flexible shifted CMRH method, the solution to
\begin{equation}
{\bm z}^{(\sigma_i)}_j = \arg\min\limits_{{\bm z}\in\mathbb{C}^{j}}\|((\check{H}_j - \underline{I}_j)
\Gamma^{(i)}_j + \underline{I}_j){\bm z} - \alpha{\bm e}_1\|_2,\quad {\bm y}^{(\sigma_i)}_j = Z^{(\sigma_i)}_j
{\bm z}^{(\sigma_i)}_j
\label{eq4.5}
\end{equation}
yields approximate solutions to Eq.~(\ref{eq1.1}) in the search spaces $Z^{(\sigma_i)}_j\in \mathbb{C}^{2n\times
j}$ that minimize the 2-norm of the \textit{quasi-residual} of the $i$-th shifted linear system, cf. Section \ref{sec3}
and \cite[Section 3.4]{KZCG,MBMBG}.  In Eq.~(\ref{eq4.5}), the Hessenberg matrix $\check{H}_j$ corresponds to the seed
system, and $\Gamma^{(i)}_j$ is related to Eq.~(\ref{eq4.4}). Note that the shifted Hessenberg procedure yields collinear residuals by default \cite{XMGTZHG}. We summarize the new nested Krylov subspace (dubbed \texttt{Hessen-FCMRH})
method for solving shifted linear systems in Algorithm \ref{alg3}.
\begin{algorithm}[!htpb]
\caption{\texttt{The FCMRH} method with pivoting and \texttt{msHessen} preconditioner}
\begin{algorithmic}[1]
\Require
the coefficient matrix $A$ (or a user-defined function that applies $A$ to a vector); the right-hand side vector ${\bm b}$; the set of shifts $\{\sigma_i\}_{i=1,\ldots,t}$; the number of inner iterations, $m_i$; the maximum number of outer iterations, $m_o$.
\Ensure
the set of solution vectors ${\bm x}^{(i)}$ of the sequence of multi-shifted linear systems.
  \State Choose the initial guess ${\bm x}_0 = {\bm x}^{(\sigma_i)}_0 = {\bm 0}$, then ${\bm r}_0 = {\bm b}$
  \State Set ${\bm q} = [1,2,\ldots n]^T$ and determine $j_0$ such that $|({\bm r}_0)_{j_0}|
         = \|{\bm r}_0\|_{\infty}$
  \State Set $\beta = ({\bm r}_0)_{j_0},\ {\bm l}_1 = {\bm r}_0/\beta$ and $({\bm q})_1 \leftrightarrow
         ({\bm q})_{j_0}$, where $\leftrightarrow$ is used to swap contents.
  \For{$j = 1,2,\ldots,m_i$,}
  \State Preconditioning: ${\bm z}^{(\sigma_i)}_j = \mathtt{msHessen}(A - \sigma_iI, {\bm l}_j,m_i)$
  \State Compute $\gamma^{(\sigma_i)}_j$ according to the similar formula in \cite[Eq. (3.9)]{MBMBG}
  \State Compute ${\bm u} = A{\bm z}^{(0)}_j$ (hint: $\sigma_i = 0$)
  \For{$k = 1,2,\ldots,j$,}
  \State $\check{h}_{k,j} = ({\bm u})_{(\bm q)_k}$
  \State ${\bm u} = {\bm u} - \check{h}_{k,j}{\bm l}_k$
  \EndFor
  \If{$j < n$ and ${\bm u}\neq {\bm 0}$}
  \State Determine $j_{0}\in \{j+1,\ldots,n\}$ such that $ |({\bm u})_{({\bm q})_{j_0}}| = \|({\bm
         u})_{({\bm q})_{j + 1}:({\bm q})_n}\|_{\infty}$
  \State $\check{h}_{j + 1,j} = ({\bm u})_{({\bm q})_{j_0}}$,\ \ ${\bm l}_{j+1} = {\bm u}/\check{
         h}_{j + 1,j};\ ({\bm q})_{j + 1} \leftrightarrow ({\bm q})_{j_0}$
  \Else
  \State $\check{h}_{j + 1,j} = 0$; Stop
  \EndIf
  \State {\small \textit{//\ \ Loop over shifted systems:}}
  \For{$i = 1,2,\ldots,t_s$}
  \State Define the matrix $Z^{(\sigma_i)}_j = [{\bm z}^{(\sigma_i)}_1,{\bm z}^{(\sigma_i)}_2,\ldots,{\bm z}^{(\sigma_i)}_j]$
  \State Set up $\check{H}_j(\sigma_i)$ according to Eq. (\ref{eq4.4})
  \State Solve ${\bm y}^{(\sigma_i)}_j = \arg\min\limits_{{\bm y}\in \mathbb{C}^{j}}\|\beta{\bm e}_1
  - \check{H}_j(\sigma_i){\bm y}\|$ with ${\bm e}_1 = [1,0,\ldots,0]^T\in \mathbb{R}^{j + 1}$
  \State ${\bm x}^{(\sigma_i)}_j = {\bm x}^{(\sigma_i)}_0 + Z^{(\sigma_i)}_j {\bm y}^{(\sigma_i)}_j$
  \EndFor
  \EndFor
\end{algorithmic}
\label{alg3}
\end{algorithm}

The proposed \texttt{Hessen-FCMRH} method is related to the FOM-FGMRES method introduced in \cite{MBMBG}.
On the other hand, according to our recent work \cite{XMGTZHG} \texttt{msHessen} can be cheaper than
\texttt{msFOM} in terms of elapsed CPU time when it is used as the inner solver. Our numerical results
presented in the next section confirm that shifted CMRH often requires less operations than shifted GMRES.
Thus, the proposed \texttt{Hessen-FCMRH} can be a cost effective nested Krylov solvers for multi-shifted
linear systems. However, note that the framework presented above may accommodate the use of other shifted
Krylov subspace methods as  inner preconditioners at Line~5 of Algorithm \ref{alg3}, for example shifted
BiCGStab($\ell$), shifted IDR($s$), shifted BiCRStab and shifted GPBiCG) could be employed. Similarly to
the FCMRH method for unshifted systems, extra memory is required to store the $Z^{(\sigma_i)}_j$'s matrices
which span the solution space for each shifted problem~\cite{YSaadA,VSDBSF}. The memory requirements increase
to about $\mathcal{O}(2ntm_o)$ for $Z^{(\sigma_i)}_{m_o}$, and each iteration of the nested multi-shifted
Krylov solver needs to evaluate $m_im_o$ matrix-vector products (the leading cost), where $m_o$ and $m_i$
are the numbers of outer and inner iterations, respectively. This extra cost is the price to pay for using
flexible preconditioning\footnote{One application of flexible preconditioning requires to evaluate
$m_i$ matrix-vector products and about $\mathcal{
O}(tm^{2}_{i})$ operations to solve $t$ Hessenberg least-square problems of size $m_i$.}~\cite{YSaadA,KZCG}, and it applies for each shift.
\section{Numerical experiments}
\label{sec5}
In this section, we show the numerical behaviour of the restarted shifted CMRH
method, shortly referred to as \texttt{sCMRH($m$)}, for solving some realistic
shifted linear systems arising in real-world engineering modelling, also compared
to the restarted shifted GMRES (\texttt{sGMRES($m$)}), restarted shifted Simpler
GMRES (\texttt{sSGMRES($m$)}) and shifted QMRIDR($s$) (\texttt{sQMRIDR($s$)})
methods. Additionally, some experiments are reported with the framework of nested
Krylov subspace solvers based on the CMRH and Hessenberg methods described in
Section~\ref{sec4}, and compared to the framework proposed in~\cite{MBMBG}.

Our experiments are performed in double precision floating point arithmetic in MATLAB
R2016a on a Windows 7 (64 bit) PC equipped with an Intel(R) Core(TM) i5-2400 CPU
running at 3.10 GHz and with 10 GB of RAM. Unless stated otherwise, the right-hand
side ${\bm b}$ is the vector with all $1$'s. The iterative solution is started from
${\bm x}^{(i)}_0 = {\bm 0} $ and is stopped at iteration $k$ when $\|{\bm r}^{(i)}_k
\|_2/\|{\bm b}\|_2 < 10^{-8}$ for $i = 1,2,\ldots,t_s$ for all linear systems, or after
at most $Max_{mvps}$ iterations. We do not compute the residuals of the additional
shifted systems because of the collinearity condition.
\subsection{General academic problems}
\label{subsec5.1}
The first set of eight linear systems are extracted from the \texttt{SuiteSparse Matrix
Collection}~\cite{TDYHT}. The main characteristics of the test problems are listed in
Table~\ref{tab1}. We use shifts $\sigma_j = -j/10000$ ($j = 1,2,3,4,5$) for problems $\Sigma_1$,
$\Sigma_4$, $\Sigma_5$, $\Sigma_6$, $\Sigma_8$, shifts $\sigma_j = -(8 + j)/200000$ for
$\Sigma_2$, shifts $\sigma_j = -(179 + j)/20000$ for $\Sigma_3$, and $\sigma_j = -j/20000$
for $\Sigma_7$. The first linear system $(A - \sigma_1 I){\bm x}^{(1)} = {\bm b}$ is the
\textit{seed} system. We set the restart value $m$ equal to $40$ (this value is used throughout Sections
\ref{subsec5.1}-\ref{sec5.4}) to ensure that all the multi-shifted linear systems converge successfully \cite{AKSTB,MBMBG,GGRG,DDRBMW},
and the maximum number of matrix-vector products $Max_{mvps}$ equal to $6000$.
\begin{table}[!htpb]
\caption{{\small Set and characteristics of the test problems used for Experiment 5.1.}}
\centering
\begin{tabular}{llrrr}
\toprule
Index   &Matrix                     & Size    &Field                         &$nnz(A)$    \\
\hline
$\Sigma_1$ &\texttt{poisson3Da}        &13,514   &Computational fluid dynamics  &352,762  \\
$\Sigma_2$ &\texttt{epb1}              &14,734   &Thermal problem               &95,053    \\
$\Sigma_3$ &\texttt{waveguide3D}       &21,036   &Electromagnetics problem      &303,468    \\
$\Sigma_4$ &\texttt{kim1}              &38,415   &2D/3D problem                 &933,195            \\
$\Sigma_5$ &\texttt{poisson3Db}        &85,623   &Computational fluid dynamics  &2,374,949          \\
$\Sigma_6$ &\texttt{vfem}              &93,476   &Electromagnetics problem      &1,434,636         \\
$\Sigma_7$ &\texttt{matrix-new\_3}     &125,329  &Semiconductor device problem  &893,984   \\
$\Sigma_8$ &\texttt{FEM\_3D\_thermal2} &147,900  &Thermal problem               &3,489,300          \\
\bottomrule
\end{tabular}
\label{tab1}
\end{table}
\begin{table}[t]\small\tabcolsep=6pt
\begin{center}
\caption{Convergence results with different shifted Krylov subspace solvers for Experiment 5.1,
using $m = 40$. Symbol $\ddag$ (or ${\cdot}^{+}$) means convergence failure (in terms of true residuals).}
\vspace{1mm}
\begin{tabular}{crrrrrrrrrr}
\hline &\multicolumn{2}{c}{\texttt{sGMRES($m$)}} &\multicolumn{2}{c}{\texttt{sCMRH($m$)}} &\multicolumn{2}
{c}{\texttt{sQMRIDR(1)}}&\multicolumn{2}{c}{\texttt{sQMRIDR(2)}}&\multicolumn{2}{c}{\texttt{sSGMRES($m$)}}
\\
[-2pt]\cmidrule(l{0.7em}r{0.7em}){2-3} \cmidrule(l{0.7em}r{0.6em}){4-5}\cmidrule(l{0.7em}r{0.7em}){6-7}
\cmidrule(l{0.7em}r{0.7em}){8-9}\cmidrule(l{0.7em}r{0.7em}){10-11} \\[-11pt]
\texttt{Index}    &\texttt{MVPs} &\texttt{CPU} & $\texttt{MVPs}$ &\texttt{CPU} &\texttt{MVPs} &
$\texttt{CPU}$&\texttt{MVPs} &\texttt{CPU} &$\texttt{MVPs}$ &$\texttt{CPU}$   \\
\hline
$\Sigma_1$    &320    &0.665    &360   &0.358  &250   &0.464  &225    &0.510  &240$^{+}$  &0.712  \\
$\Sigma_2$    &1160   &1.339    &1769  &0.963  &1367  &1.694  &989    &1.559  &1050$^{+}$ &2.444  \\
$\Sigma_3$    &480    &2.644    &560   &2.065  &437   &2.640  &392    &3.275  &\ddag      &\ddag  \\
$\Sigma_4$    &2000   &20.390   &1640  &11.268 &\ddag &\ddag  &2231   &36.052 &\ddag      &\ddag  \\
$\Sigma_5$    &600    &6.741    &680   &6.211  &561   &9.309  &427    &9.527  &484$^{+}$  &9.524  \\
$\Sigma_6$    &240    &4.461    &280   &3.367  &413   &10.241 &224    &9.189  &\ddag      &\ddag  \\
$\Sigma_7$    &120    &0.789    &120   &0.476  &287   &3.423  &152    &3.215  &\ddag      &\ddag  \\
$\Sigma_8$    &520    &12.124   &560   &6.764  &793   &17.560 &503    &16.932 &\ddag      &\ddag  \\
\hline
\end{tabular}
\label{tab2}
\end{center}
\end{table}

In Table \ref{tab2} we show number of matrix-vector products (abbreviated as \texttt{MVPs}) and elapsed CPU solution time in seconds (\texttt{CPU}) required by \texttt{sCMRH($m$)}, \texttt{sGMRES($m$)}, \texttt{sQMRIDR($s$)}$(s = 1,2)$ and \texttt{sSGMRES($m$)} to converge to prescribed accuracy. All methods converge within the maximum number of iterations, except \texttt{sQMRIDR($1$)} for Problem $\Pi_4$ and \texttt{sSGMRES($m$)} for Problem $\Pi_{\ell},~\ell = 3,4,6,7,8$. In our runs, \texttt{sCMRH($m$)} is the fastest solver although it often needs more \texttt{MVPs}, confirming the complexity analysis presented in Section~\ref{sec3.2}. The timing performance of \texttt{sQMRIDR($s$)} are penalized by  the extra inner products and vector updates required at each iteration; note, however, that \texttt{sQMRIDR($2$)} generally requires less matrix-vector products to converge. On this set of problems, \texttt{sSGMRES($m$)} is less robust that the other solvers.

\begin{figure}[H]
\centering
\includegraphics[width=3.12in,height=2.95in]{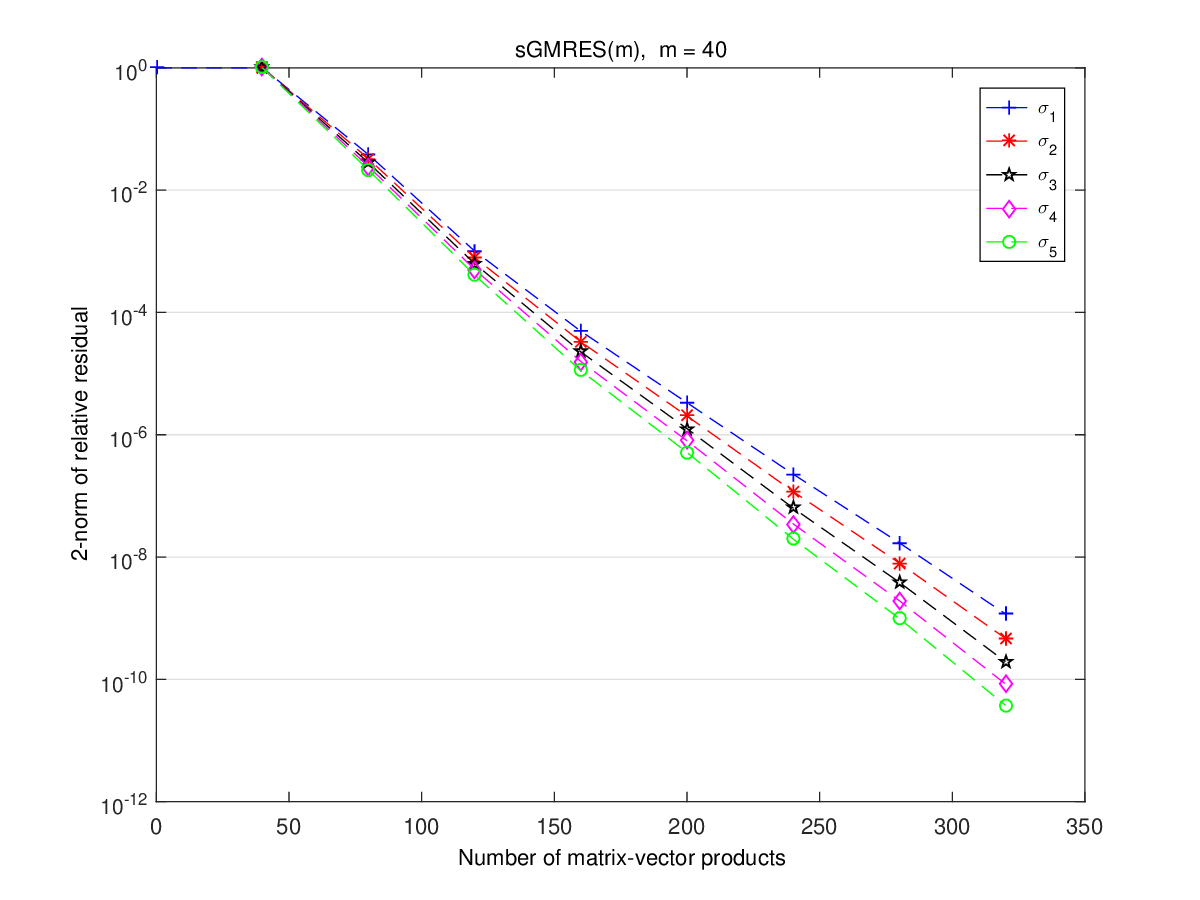}
\includegraphics[width=3.12in,height=2.95in]{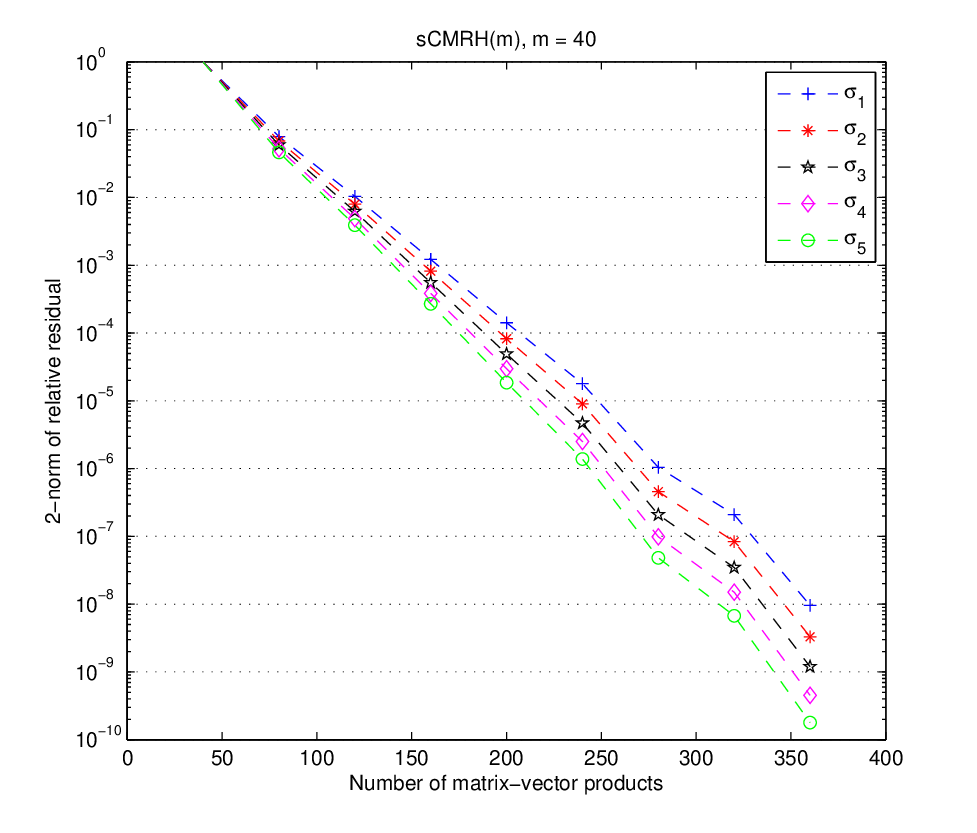}
\includegraphics[width=3.12in,height=2.95in]{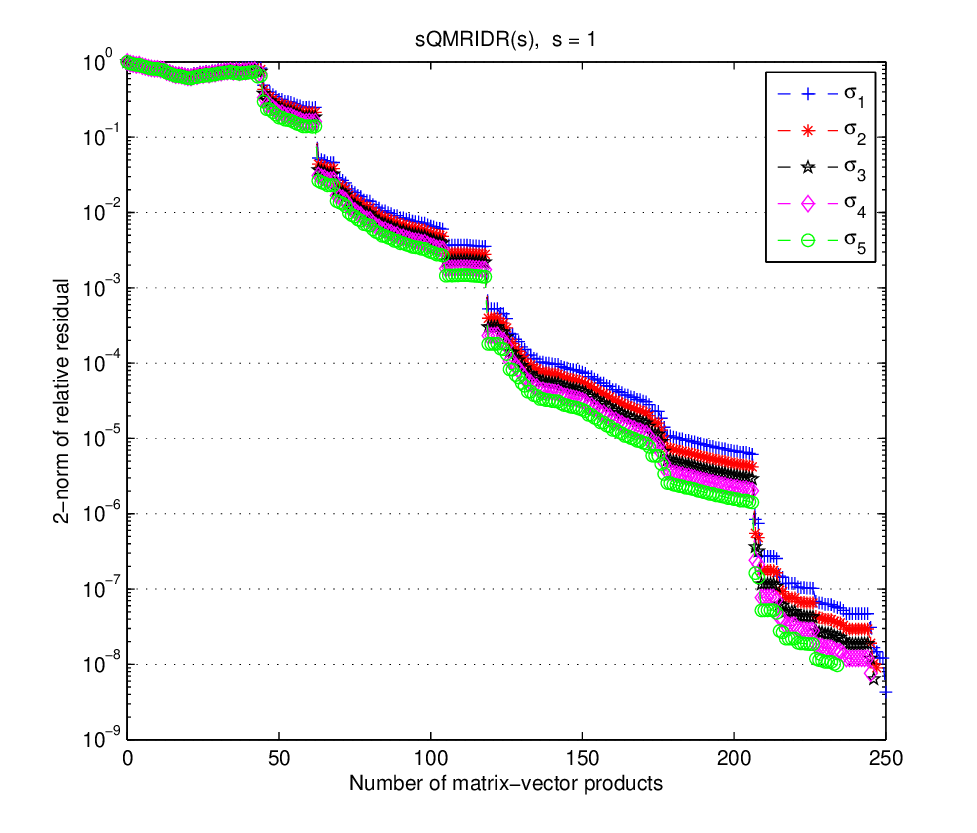}
\includegraphics[width=3.12in,height=2.95in]{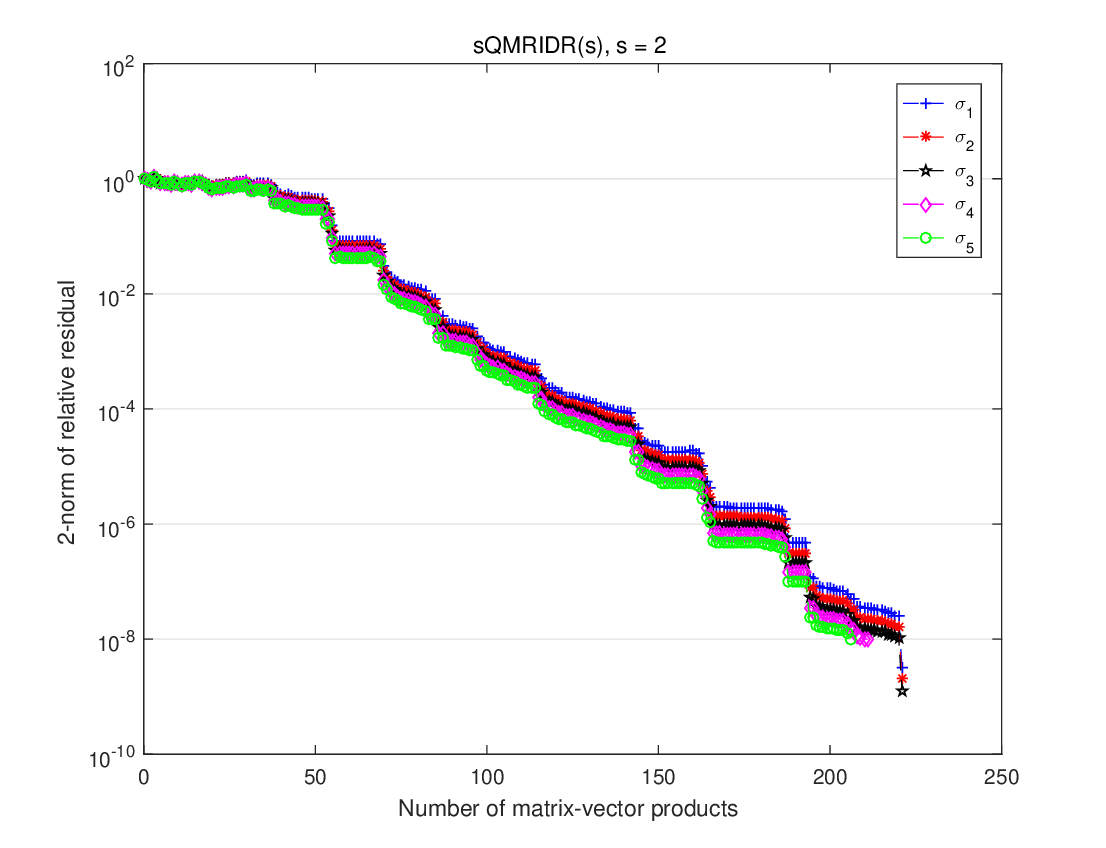}
\caption{{\small Relative residual histories of different iterative solvers for
test problem $\Sigma_1$ in Table \ref{tab2}.}}
\label{fig1}
\end{figure}

The relative residual histories plotted in Fig.~\ref{fig1} for different shifts confirm the conclusion reported in~\cite{AFUG} that \texttt{sGMRES($m$)} converges more rapidly on the {\it add} systems than on the seed system, and that its convergence history is comparable  to \texttt{sCMRH($m$)}. Overall, the \texttt{sCMRH($m$)} method is very competitive on this set of matrices.
\subsection{Quantum chromodynamics applications}
\label{sec:5.2}
Next, we report on experiments on seven $49,152\times49,152$ complex shifted linear systems arising
from the discretization of the Dirac operator in quark simulations at different physical temperatures
in lattice QCD applications. The seven test problems are named as \texttt{conf5\_4-8x8-05}, \texttt{conf5\_4-8x8-10},
\texttt{conf5\_4-8x8-15}, \texttt{conf5\_4-8x8-20}, \texttt{conf6\_0-8x8-20}, \texttt{conf6\_0-8x8-30},
\texttt{conf6\_0-8x8-80} in the \texttt{SuiteSparse Matrix Collection}~\cite{TDYHT}. Hereafter, they will be denoted as problems $\Pi_i$ for $i = 1,2,\ldots,7$. We select shifts equals to $\sigma_j\in \mathcal{I} = -\{.001, .002, .003, .04, .05, .06, .07\}$ and we use $(A - \sigma_1 I){\bm x}^{(1)} = {\bm b}$ as the \textit{seed} system. The maximum number of matrix-vector products is set equal to $5000$.
\begin{table}[!htpb]\small\tabcolsep=6pt
\begin{center}
\caption{Convergence results with different shifted Krylov subspace solvers for Experiment~5.2, using $m = 40$.}
\vspace{1mm}
\begin{tabular}{crrrrrrrrrr}
\hline &\multicolumn{2}{c}{\texttt{sGMRES($m$)}} &\multicolumn{2}{c}{\texttt{sCMRH($m$)}} &\multicolumn{2}
{c}{\texttt{sQMRIDR(1)}}&\multicolumn{2}{c}{\texttt{sQMRIDR(2)}}&\multicolumn{2}{c}{\texttt{sAd-SGMRES($m$)}}
\\
[-2pt]\cmidrule(l{0.7em}r{0.7em}){2-3} \cmidrule(l{0.7em}r{0.6em}){4-5}\cmidrule(l{0.7em}r{0.7em}){6-7}
\cmidrule(l{0.7em}r{0.7em}){8-9}\cmidrule(l{0.7em}r{0.7em}){10-11} \\[-11pt]
\texttt{Index}    &\texttt{MVPs} &\texttt{CPU} & $\texttt{MVPs}$ &\texttt{CPU} &\texttt{MVPs} &
$\texttt{CPU}$&\texttt{MVPs} &\texttt{CPU} &$\texttt{MVPs}$ &$\texttt{CPU}$   \\
\hline
$\Pi_1$    &1120  &15.007   &1480   &14.974 &845  &15.557 &732  &20.549 &3459$^{+}$ &82.897  \\
$\Pi_2$    &880   &11.388   &1120   &11.236 &746  &13.835 &675  &19.456 &\ddag      &\ddag  \\
$\Pi_3$    &720   & 9.173   &880    & 8.815 &709  &13.345 &633  &19.018 &\ddag      &\ddag  \\
$\Pi_4$    &640   & 8.265   &800    & 7.986 &757  &13.973 &687  &20.283 &\ddag      &\ddag  \\
$\Pi_5$    &640   & 8.255   &840    & 8.279 &589  &11.473 &477  &14.797 &\ddag      &\ddag  \\
$\Pi_6$    &840   &10.762   &1160   &11.712 &539  &10.257 &456  &13.468 &\ddag      &\ddag  \\
$\Pi_7$    &680   & 8.656   &840    & 8.321 &582  &11.171 &477  &14.109 &1743$^{+}$ &44.703  \\
\hline
\end{tabular}
\label{tab3}
\end{center}
\end{table}

In Table \ref{tab3}, we report on elapsed CPU time and \texttt{MVPs} required by different iterative methods to reduce the initial residuals by eight orders of magnitude. Except for Problems $\Pi_5$ and $\Pi_6$, \texttt{sCMRH($m$)}  is more cost-effective that \texttt{sGMRES($m$)}. In general \texttt{sQMRIDR($s$)} requires less iterations but more CPU time to converge. One exception is Problem $\Pi_6$, where \texttt{sQMRIDR(1)} is the fastest solver. In our experiments, \texttt{sAd-SGMRES($m$)} is not a competitive choice in terms of both \texttt{MVPs} and elapsed CPU time.

\begin{figure}[!htpb]
\centering
\includegraphics[width=3.12in,height=2.95in]{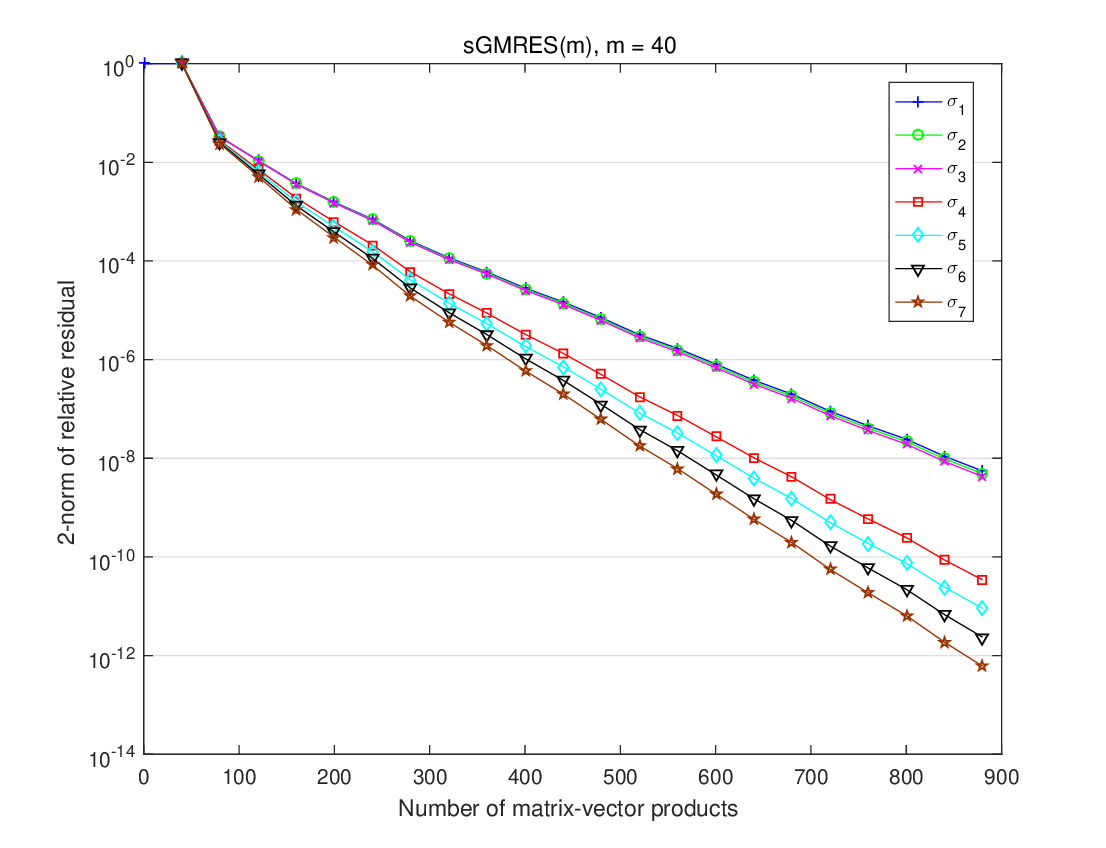}
\includegraphics[width=3.12in,height=2.95in]{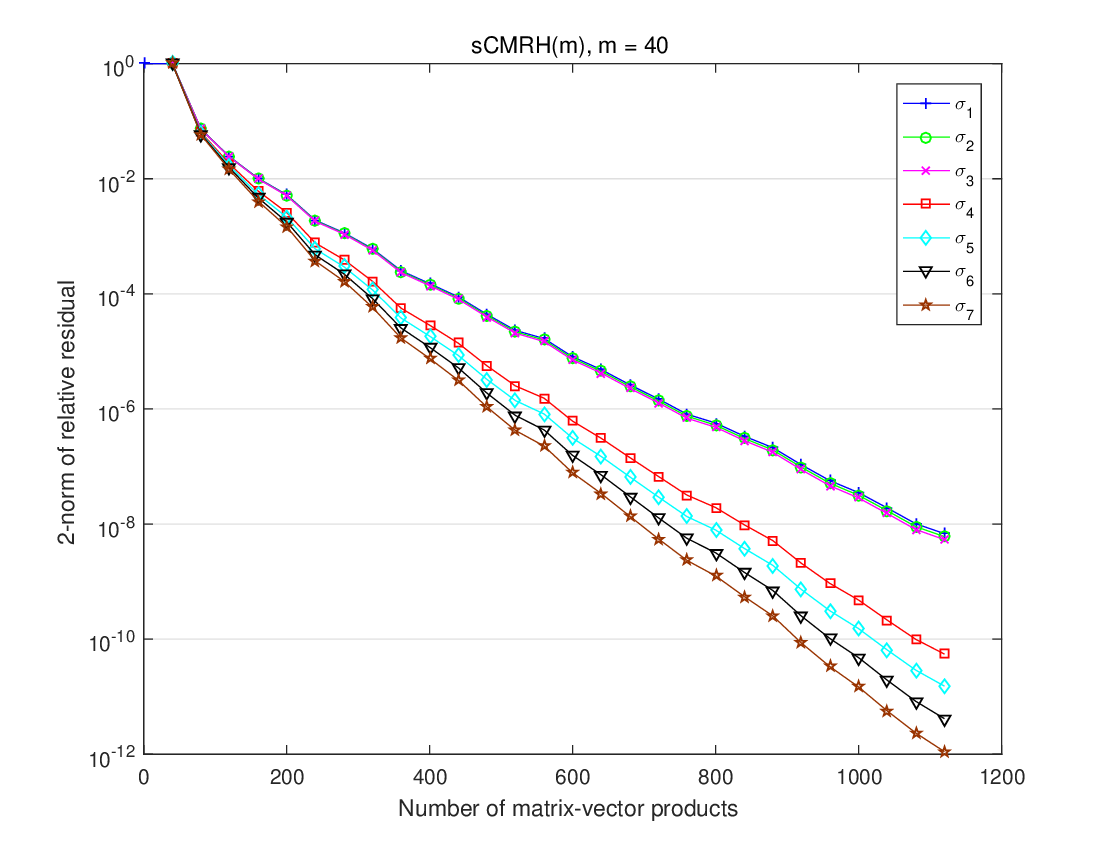}
\includegraphics[width=3.12in,height=2.95in]{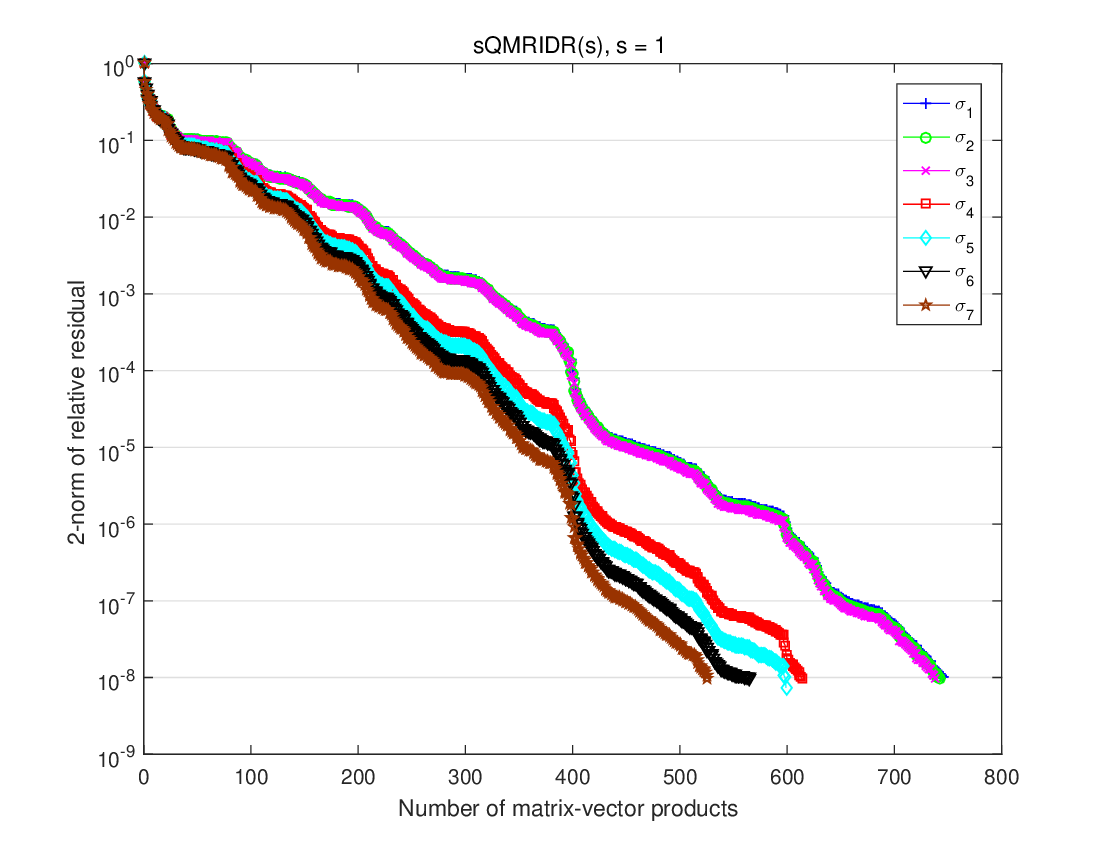}
\includegraphics[width=3.12in,height=2.95in]{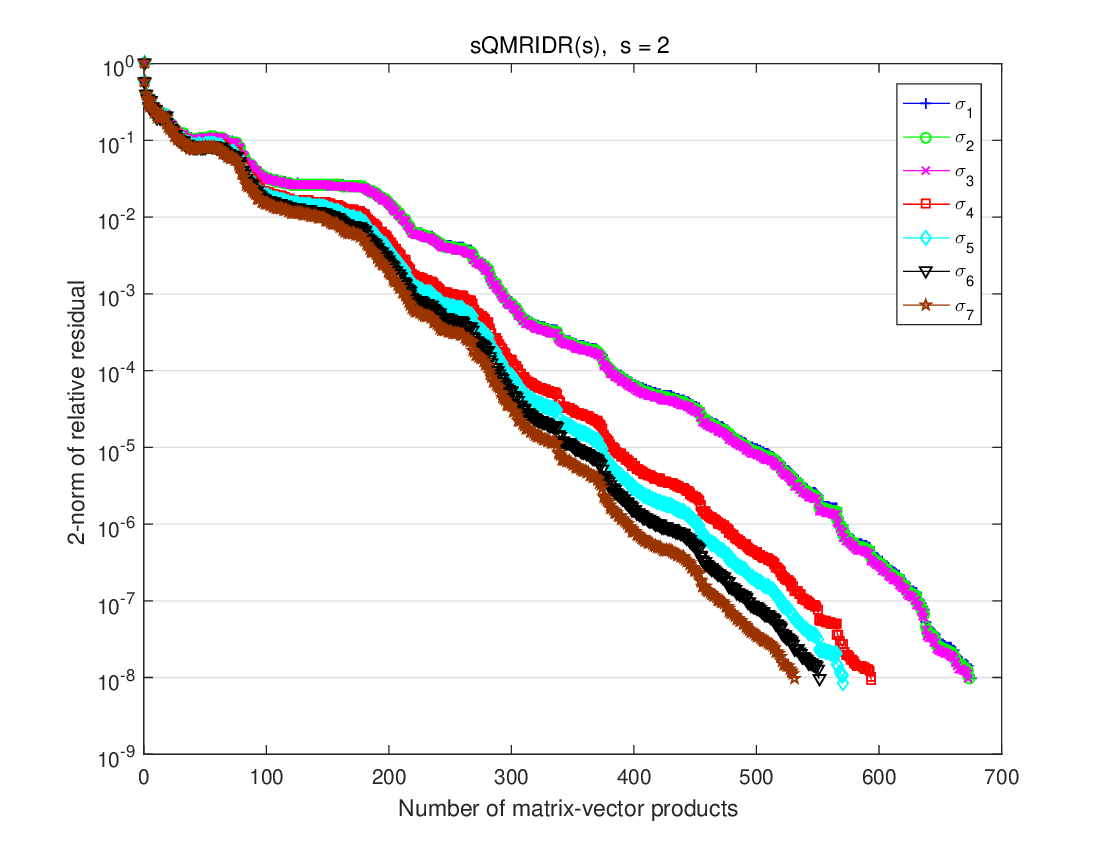}
\caption{{\small Relative residual histories of different iterative solvers for
test problem $\Pi_2$ in Table \ref{tab3}.}}
\label{fig2}
\end{figure}

The convergence histories and the final accuracies of \texttt{sGMRES($m$)} and \texttt{sCMRH($m$)} are once again comparable based on the relative residual histories plotted in Fig. \ref{fig2}. The approximate solutions computed by \texttt{sCMRH($m$)} are slightly more accurate than for the other solvers, whereas the final residual norms of \texttt{sGMRES($m$)} are smaller than those of \texttt{sQMRIDR($s$)}. We conclude that \texttt{sCMRH($m$)} can be an interesting alternative to other shifted Krylov subspace methods for this set of problems.
\subsection{Time fractional differential equations}
\label{sec5.4}
Fractional differential equations (FDEs) are widespread mathematical models in the study of physical, biological, geological and financial systems, to name only a few fields. Due to  the increasing dimension and complexity of these systems, in recent years considerable attention has been devoted to the development of efficient numerical methods for the approximate solution of FDEs in many  areas, see e.g. \cite{MCCES,IPFD}. Here we consider a benchmark problem coming from the 3D time-fractional convection-diffusion-reaction equation defined as
\begin{equation}
\begin{cases}
\frac{\partial^{\gamma} u}{\partial t^{\gamma}} = \epsilon \triangle u - \vec{\beta}\cdot
\nabla u + r u, & (x,y,z)\in \Omega = (0,1)^3,\ t \in [0,T],\\
u(x,y,z,t) = 0, & (x,y,z)\in \partial\Omega,\ t \in [0,T], \\
u(x,y,z,0) = x(1 - x)y(1-y)z(1 - z),& (x,y,z)\in \bar{\Omega}.
\end{cases}
\label{eq5.1}
\end{equation}
Problem~(\ref{eq5.1}) is a modification of the third example presented in Ref. \cite{MBGGLG}, with diffusion
constant $\epsilon = 1$, reaction constant $r = 400$ and convection parameter $\vec{\beta} = (0/\sqrt{5},
250/\sqrt{5}, 500/\sqrt{5})^T$. Upon the finite difference discretization of Eqs.~(\ref{eq5.1}) on an uniformly
spaced domain using naturally ordered grid points, a system of FDEs of the form
\begin{equation}
\frac{d^{\gamma} {\bm u}}{d t^{\gamma}} = A{\bm u}(t),\qquad {\bm u}(0)={\bm u}_0.
\label{eq5.2}
\end{equation}
is obtained, where ${\bm u}$ denotes the vector of unknown approximate solutions at the grid points. Using a grid
size $h = 0.025$, the order of the matrix $A$ is about 60,000. It is well-known \cite{LNTJAC,JACWL,RGAFMP} that, for $0 < \gamma < 1$, the true solution of this problem can be expressed as
\begin{equation}
{\bm u}(t) = e_{\gamma,1}(t;A){\bm u}_0,\quad\ \mathrm{and}\quad\ e_{\gamma,1}(t;A) =
t^{1-1}E_{\gamma,1}(t^{\gamma}A)= E_{\gamma,1}(t^{\gamma}A),
\label{eq5.3}
\end{equation}
where $E_{\gamma,1}(z)$ is the Mittag-Leffler (ML) function \cite{IPFD,RGAFMP}
\begin{equation*}
E_{\gamma,1}(z) := \sum^{\infty}_{k=0}\frac{z^k}{\Gamma(\gamma k + 1)},\quad\
\gamma > 0,\ \ z\in\mathbb{C}.
\end{equation*}
In light of Eq.~(\ref{eq5.3}), the numerical solution ${\bm u}(t)$ can be computed as the product of the matrix ML function $E_{\gamma,1}(t^{\gamma}A)$ times ${\bm u}_0 $. This operation  accounts for the major computational cost of the solution of Eq.~(\ref{eq5.2}). Recently, the numerical evaluation of the action of the matrix function on a vector, namely $E_{\gamma,1}(t^{\gamma}A){\bm u}_0$, is receiving much consideration as shown by the spread of literature on this topic~\cite{LNTJAC,JACWL,RGAFMP}.
One approach is based on the Carath\'{e}odory-Fej\'{e}r approximation of $E_{\gamma,1}(t^{\gamma}A){\bm u}_0$~\cite{LNTJAC2}, that has the following representation
\begin{equation}
E_{\gamma,1}(A){\bm u}_0 = f_{\nu}(A){\bm u}_0 = \sum^{\nu}_{j = 1}w_j(z_j I - A)^{-1}{\bm u}_0,
\label{eq5.4}
\end{equation}
where $w_j$ and $z_j$ are quadrature weights and nodes, respectively. Implementing Eq.~(\ref{eq5.4})
requires the solution of a sequence of shifted linear systems of the form $(-A + z_j I) {\bm x}^{(j)}={\bm u}_0,\ z_j\in \mathbb{C}$.
\begin{table}[t]\small\tabcolsep=6pt
\begin{center}
\caption{Convergence results with different shifted Krylov subspace solvers for Experiment 5.3, using $m = 40$. Symbol $\ddag$ means convergence failure.}
\vspace{2mm}
\begin{tabular}{rcccccrcrcr}
\hline &\multicolumn{2}{c}{\texttt{sGMRES($m$)}} &\multicolumn{2}{c}{\texttt{sCMRH($m$)}}
&\multicolumn{2}{c}{\texttt{sQMRIDR(1)}}&\multicolumn{2}{c}{\texttt{sQMRIDR(2)}}&\multicolumn{2}
{c}{\texttt{sAd-SGMRES($m$)}}
\\
[-2pt]\cmidrule(l{0.7em}r{0.7em}){2-3} \cmidrule(l{0.7em}r{0.6em}){4-5}\cmidrule(l{0.7em}r{0.7em}){6-7}
\cmidrule(l{0.7em}r{0.7em}){8-9}\cmidrule(l{0.7em}r{0.7em}){10-11} \\[-11pt]
$(\gamma,\nu)$ &\texttt{MVPs} &\texttt{CPU} & $\texttt{MVPs}$ &\texttt{CPU} &\texttt{MVPs} &
$\texttt{CPU}$&\texttt{MVPs} &\texttt{CPU} &$\texttt{MVPs}$ &$\texttt{CPU}$   \\
\hline
$(0.2,6)$  &800  &3.421   &1000   &2.990   &299  &6.225  &182   &6.081  &\ddag  & \ddag  \\
$(0.4,8)$  &800  &3.553   &1000   &3.051   &299  &8.061  &182   &7.809  &\ddag  & \ddag \\
$(0.6,10)$ &800  &3.560   &1000   &3.130   &299  &9.651  &182   &9.542  &\ddag  & \ddag \\
$(0.8,10)$ &800  &3.547   &1000   &3.136   &299  &9.755  &182   &9.633  &\ddag  & \ddag \\
$(0.9,12)$ &800  &3.589   &1000   &3.204   &299  &11.302 &185   &11.317 &\ddag  & \ddag \\
\hline
\end{tabular}
\label{tab6}
\end{center}
\end{table}
\begin{table}[!htbp]\small\tabcolsep=6pt
\begin{center}
\caption{Convergence results with different shifted Krylov solvers for Experiment 5.3, using $m = 40$, $\gamma = 0.8$ and $\nu
= 10$.}
\vspace{1mm}
\begin{tabular}{lcccrcrcrcr}
\hline &\multicolumn{2}{c}{\texttt{sGMRES($m$)}} &\multicolumn{2}{c}{\texttt{sCMRH($m$)}} &\multicolumn{2}
{c}{\texttt{sQMRIDR(1)}}&\multicolumn{2}{c}{\texttt{sQMRIDR(2)}}&\multicolumn{2}{c}{\texttt{sAd-SGMRES($m$)}}
\\
[-2pt]\cmidrule(l{0.7em}r{0.7em}){2-3} \cmidrule(l{0.7em}r{0.6em}){4-5}\cmidrule(l{0.7em}r{0.7em}){6-7}
\cmidrule(l{0.7em}r{0.7em}){8-9}\cmidrule(l{0.7em}r{0.7em}){10-11} \\[-11pt]
$(h,~r)$ &\texttt{MVPs} &\texttt{CPU} & $\texttt{MVPs}$ &\texttt{CPU} &\texttt{MVPs} &
$\texttt{CPU}$&\texttt{MVPs} &\texttt{CPU} &$\texttt{MVPs}$ &$\texttt{CPU}$   \\
\hline
$(\frac{1}{50},400)$  &$\ddag$ &$\ddag$ &1360 &8.41   &283 &22.52  &200 &24.44  &2653$^{+}$ &197.72  \\
$(\frac{1}{50},500)$  &$\ddag$ &$\ddag$ &1400 &8.63   &293 &23.58  &215 &25.94  &\ddag      &\ddag   \\
$(\frac{1}{80},400)$  &$\ddag$ &$\ddag$ &1640 &61.57  &358 &153.55 &275 &162.09 &\ddag      &\ddag   \\
$(\frac{1}{80},500)$  &$\ddag$ &$\ddag$ &3000 &114.90 &359 &155.35 &278 &162.61 &\ddag      &\ddag   \\
$(\frac{1}{100},400)$ &$\ddag$ &$\ddag$ &1680 &138.81 &397 &355.76 &329 &391.61 &\ddag      &\ddag   \\
$(\frac{1}{100},500)$ &$\ddag$ &$\ddag$ &2680 &228.47 &402 &359.26 &333 &395.71 &\ddag      &\ddag   \\
\hline
\end{tabular}
\label{tab6xx}
\end{center}
\end{table}

In Table~\ref{tab6} we compare different shifted Krylov subspace methods for solving five groups of real shifted linear systems with same seed system $-A{\bm x} = {\bm u}_0$ but a different number of shifts arising from this application. For this reason, the number of matrix-vector products required by different solvers to converge changes only slightly in Table~\ref{tab6}. In our runs, \texttt{sCMRH($m$)} exhibits the fastest  convergence in terms of elapsed CPU time, although it requires more \texttt{MVPs}. As in previous experiments, the performance of \texttt{sQMRIDR($s$)} are penalized by the $s + 2$, possibly complex, $n$-length extra vectors that need to be stored and updated for each new shift. Unfortunately, the \texttt{sAd-SGMRES($m$)} method fails to converge in our runs.

For a more comprehensive performance evaluation of different algorithms, in Table \ref{tab6xx} we report on elapsed CPU time and \texttt{MVPs} required to solve the six shifted linear systems with different grid size and reaction parameters.  Large values of convection coefficients $\vec{\beta}$ and reaction parameter $r$ in Eq. (\ref{eq4.1}) always result in highly
nonsymmetric, ill-conditioned finite difference matrix $A$ \cite{MBGGLG}. In these experiment we set $Max_{mvps}= 6000$. We can see from Table~\ref{tab6xx} that in our runs only \texttt{sCMRH($m$)} and \texttt{sQMRIDR($s$)} are able to solve all the shifted linear systems. The reason might be that the \texttt{sGMRES($m$)} and \texttt{sAd-SGMRES($m$)} methods are more sensitive than the other two shifted iterative solvers to the high nonsymmetry and indefiniteness of the coefficient matrix. In addition, although \texttt{sQMRIDR($s$)} always requires less number of \texttt{MVPs} than the other two shifted iterative solvers, it is still more expensive in terms of elapsed CPU time due to the extra vector operations. Based on these results, we conclude that the proposed \texttt{sCMRH($m$)} method is an efficient algorithm to solve shifted linear systems arising in the discretization of FDEs.
\subsection{Experiments on flexible preconditioning for shifted systems}
\label{sec5x}
We illustrate the performance of the methods presented in Section~\ref{sec4}, namely \texttt{Hessen-FCMRH},
\texttt{Hessen-FGMRES} and \texttt{FOM-FCMRH} against \texttt{FOM-FGMRES} \cite{MBMBG} for solving some ill-conditioned shifted linear systems arising from realistic difficult problems, such as evaluating the action of the matrix function on a vector with large matrix norm \cite{LNTJAC2,MCPKAO} and the QCD simulation. In these runs, we select $A{\bm x}
= {\bm b}$ as the \textit{seed} system. Timings can be different from the experiments of the previous two sections as we use the norm of the \underline{\textit{true residual vector}} (which can help us to investigate how inner iterations actually affect the outer iterations) to monitor the (complete) convergence of the outer method (i.e., \texttt{FCMRH} and \texttt{FGMRES}).

In \texttt{Group \uppercase\expandafter{\romannumeral1}}, we solve some shifted linear systems arising from the computation of the Carath\'{e}odory-Fej\'{e}r approximation $\exp(\tau L){\bm u}_0$  in reactive transport simulations through heterogeneous porous media modelled by the advection-diffusion-reaction equation defined on a two-dimensional $[0,2]^2$ square domain. The experimental setting is the same as the one proposed in \cite[Example 2]{MCPKAO}, except for different values of the shift $\tau$. We use MATLAB codes available at
\url{https://numerical-analysis.uibk.ac.at/exp-int-software} to generate matrices $L$ of different sizes equal to $9,801\times 9,801$, $14,161 \times 14,161$, and $19,321\times 19,321$, corresponding to grid sizes $h = 1/100,1/120,1/140$. In the matrix ML function evaluation, we use $\nu = 14$ poles. Numerical experiments with various nested Krylov subspace solvers are shown in Table~\ref{tab7}.

\begin{table}[!htbp]\small\tabcolsep=6pt
\begin{center}
\caption{Number of inner and outer iterations and \texttt{CPU} solution time for experiments with different nested Krylov methods to solve \texttt{Group \uppercase\expandafter{\romannumeral1}} problems.}
\vspace{1mm}
\begin{tabular}{crcccccrcrrr}
\hline & &\multicolumn{2}{c}{\texttt{Hessen-FCMRH}} &\multicolumn{2}{c}{\texttt{Hessen-FGMRES}} &\multicolumn{2}
{c}{\texttt{FOM-FCMRH}}&\multicolumn{2}{c}{\texttt{FOM-FGMRES}}\\
[-2pt]\cmidrule(l{0.7em}r{0.7em}){3-4} \cmidrule(l{0.7em}r{0.6em}){5-6}\cmidrule(l{0.7em}r{0.7em}){7-8}
\cmidrule(l{0.7em}r{0.7em}){9-10} \\[-11pt]
$(h,\tau)$& \texttt{IT\_in} &\texttt{IT\_out} &\texttt{CPU} & $\texttt{IT\_out}$ &\texttt{CPU} &\texttt{IT\_out} &
$\texttt{CPU}$&\texttt{IT\_out} &\texttt{CPU}   \\
\hline
$(\frac{1}{100},0.04)$  &110  &2  &0.364   &3  &0.543  &4  &1.147  &4   &1.219   \\
                        &100  &2  &0.316   &4  &0.667  &6  &1.520  &10  &2.613   \\
                        & 90  &2  &0.291   &3  &0.445  &27 &6.570  &15  &3.163   \\
$(\frac{1}{120},0.04)$  &160  &5  &2.208   &4  &1.789  &5  &3.347  &6   &4.398   \\
                        &150  &6  &2.341   &6  &2.422  &29 &20.155 &8   &4.919   \\
                        &140  &21 &8.023   &16 &5.912  &22 &13.721 &12  &6.882   \\
$(\frac{1}{140},0.03)$  &160  &4  &2.127   &5  &2.828  &4  &3.433  &3   &2.641   \\
                        &150  &6  &3.027   &6  &3.190  &4  &3.204  &4   &3.266   \\
                        &140  &6  &2.913   &8  &3.863  &7  &4.932  &10  &7.126   \\
\hline
\end{tabular}
\label{tab7}
\end{center}
\end{table}

In \texttt{Group \uppercase\expandafter{\romannumeral2}}, we use the same test problems of Section \ref{sec5.4} with the following setup values: $\gamma = 0.9$ and $\vec{\beta} = (0/\sqrt{5}, 250/\sqrt{5}, 500/\sqrt{5})^T$ with $h = 0.025$ and $h = 0.02$ for test problems $\Xi_1$ and $\Xi_2$. The sizes of these two systems are $59,319\times 59,319$ and $117,649\times 117,649$, respectively. On the other hand, we set $\gamma = 0.8$ and $\vec{\beta} = (500/\sqrt{5}, 250/\sqrt{5}, 500/\sqrt{5})^T$ with $h = 0.02$ and $h = 0.0125$ for test problems $\Xi_3$ and $\Xi_4$. The sizes of these two test problems are $117,649\times 117,649$ and $493,039\times 493,039$, respectively. The other settings are set equal to the values used in Section~5.3. The results are illustrated in Table~\ref{tab8}.

\begin{table}[t]\small\tabcolsep=6pt
\begin{center}
\caption{Number of inner and outer iterations and \texttt{CPU} solution time for experiments with different nested Krylov methods to solve \texttt{Group \uppercase\expandafter{\romannumeral2}} problems.}
\vspace{1mm}
\begin{tabular}{rrcrcrcrcrcr}
\hline & &\multicolumn{2}{c}{\texttt{Hessen-FCMRH}} &\multicolumn{2}{c}{\texttt{Hessen-FGMRES}} &\multicolumn{2}
{c}{\texttt{FOM-FCMRH}}&\multicolumn{2}{c}{\texttt{FOM-FGMRES}}\\
[-2pt]\cmidrule(l{0.7em}r{0.7em}){3-4} \cmidrule(l{0.7em}r{0.6em}){5-6}\cmidrule(l{0.7em}r{0.7em}){7-8}
\cmidrule(l{0.7em}r{0.7em}){9-10} \\[-11pt]
\texttt{Index}& \texttt{IT\_in} &\texttt{IT\_out} &\texttt{CPU} & $\texttt{IT\_out}$ &\texttt{CPU} &\texttt{IT\_out} &
$\texttt{CPU}$&\texttt{IT\_out} &\texttt{CPU}   \\
\hline
$\Xi_1$  &80   &2  &1.281   &3    &1.989   &3       &2.626   &3   &2.679    \\
         &70   &4  &2.294   &5    &2.938   &4       &2.998   &4   &3.044    \\
         &60   &18 &11.442  &7    &3.766   &6       &3.826   &5   &3.030    \\
$\Xi_2$  &90   &3  &4.344   &3    &4.368   &4       &7.843   &4   &7.847   \\
         &80   &5  &6.639   &5    &6.658   &6       &10.633  &4   &6.947   \\
         &70   &11 &13.958  &8    &9.778   &26      &48.284  &7   &10.856   \\
$\Xi_3$  &100  &2  &3.123   &3    &4.732   &5       &11.094  &4   &8.739   \\
         &90   &5  &7.106   &5    &7.187   &16      &33.751  &8   &15.893   \\
         &80   &7  &8.948   &9    &12.011  &12      &21.645  &7   &11.992   \\
$\Xi_4$  &140  &2  &29.866  &4    &60.604  &13      &559.623 &6   &252.897  \\
         &130  &6  &82.819  &7    &98.412  &18      &700.114 &4   &149.293  \\
         &120  &4  &49.197  &9    &115.089 &28      &965.549 &4   &129.557  \\
\hline
\end{tabular}
\label{tab8}
\end{center}
\end{table}

In Table \ref{tab7}, we show elapsed CPU time and  number of inner and outer iterations (denoted as \texttt{IT\_in}
and \texttt{IT\_out}, respectively) required by the four nested Krylov methods to solve the seven shifted linear systems within prescribed tolerance and maximum number of iterations. We adopt these notations throughout the subsection. Our nested solver \texttt{Hessen-FCMRH} is  the most efficient one in terms of elapsed CPU time at equal number of inner steps, with the exception of $(1/120,0.04)$ and \texttt{IT\_in} = 140. The result is supported by the conclusions presented in~\cite{XMGTZHG} and by our previous analysis. Note that not only \texttt{Hessen-FCMRH} is faster than \texttt{Hessen-FGMRES}, but also \texttt{FOM-FCMRH} is more efficient than \texttt{FOM-FGMRES} in terms of elapsed time at similar number of outer iterations. We conclude that nested iterative solvers based on the Hessenberg procedure can be computationally efficient to solve multi-shifted linear systems.

Similar conclusions are derived from the convergence results of the seven shifted linear systems of \texttt{Group \uppercase\expandafter{\romannumeral2}} illustrated in Table \ref{tab8}. Nested \texttt{Hessen-FCMRH} is a robust method in terms of solution time at equal number of inner steps, except only problem $\Xi_1$ with \texttt{IT\_in} = 60 and problem $\Xi_2$ with \texttt{IT\_in} = 70. For most problems, \texttt{Hessen-FCMRH} and \texttt{Hessen-FGMRES} converge faster than \texttt{FOM-FGMRES} and \texttt{FOM-FCMRH}. In addition, we can see that \texttt{FOM-FCMRH} is still faster than \texttt{FOM-FGMRES} in terms of elapsed time at similar number of outer iterations. Under similar conditions, \texttt{Hessen-FCMRH} is also faster than \texttt{Hessen-FGMRES} thanks to the cost-effective Hessenberg procedure. The proposed \texttt{Hessen-FCMRH} method is remarkably robust for handling shifted linear systems, in terms of the elapsed CPU time.

\begin{table}[!htbp]\small\tabcolsep=5.5pt
\begin{center}
\caption{Number of inner and outer iterations and \texttt{CPU} solution time for experiments with different nested Krylov methods to solve problem \texttt{MLMF}.}
\vspace{1mm}
\begin{tabular}{rcrcrcrcrcr}
\hline &\multicolumn{2}{c}{\texttt{IT\_in = 120}} &\multicolumn{2}{c}{\texttt{IT\_in = 110}}
&\multicolumn{2}{c}{\texttt{IT\_in = 100}}&\multicolumn{2}{c}{\texttt{IT\_in = 90}}\\
[-2pt]\cmidrule(l{0.7em}r{0.7em}){2-3} \cmidrule(l{0.7em}r{0.6em}){4-5}\cmidrule(l{0.7em}r{0.7em}){6-7}
\cmidrule(l{0.7em}r{0.7em}){8-9} \\[-11pt]
\texttt{Solver} &\texttt{IT\_out} &\texttt{CPU} & $\texttt{IT\_out}$ &\texttt{CPU} &\texttt{IT\_out} &
$\texttt{CPU}$&\texttt{IT\_out} &\texttt{CPU}   \\
\hline
\texttt{Hessen-FCMRH}      &2        &4.18   &2        &3.83   &2        &3.44   &6        &9.45   \\
\texttt{Hessen-FGMRES}     &2        &4.20   &3        &5.82   &3        &5.23   &6        &9.49   \\
\texttt{FOM-FCMRH}         &2        &5.89   &3        &7.86   &5        &11.83  &18       &43.73  \\
\texttt{FOM-FGMRES}        &2        &6.01   &3        &7.95   &5        &11.89  &8        &16.76  \\
\hline
\end{tabular}
\label{tab8xx}
\end{center}
\end{table}

Finally, we consider the last two test problems, denoted as \texttt{MLMF} and \texttt{QCDx}. \texttt{MLMF} arises from Eq. (\ref{eq4.1}) setting  $\gamma = 0.9$,
$r = 500$, $\vec{\beta} =(500/\sqrt{5},250/\sqrt{5}, 500/\sqrt{5})^T$ and $h = 0.02$. Experiments with nested iterative schemes with a  different number of inner iterations (i.e, \texttt{IT\_in}) are shown in Table \ref{tab8xx}. Our \texttt{Hessen-FCMRH} method with \texttt{IT\_out = 2} and \texttt{IT\_in =
100} exhibits the best overall performance in terms of CPU time. Nested schemes based on the Arnoldi procedure always require more time to converge, except the case of \texttt{IT\_in = 90}. Often in our runs, the larger the number of inner iterations, the smaller the number of outer steps, and vice versa for the solution time. The results show that the proposed \texttt{Hessen-FCMRH} method can be regarded as a robust choice for this problem.

\begin{table}[!htbp]\small\tabcolsep=5.5pt
\begin{center}
\caption{Number of inner and outer iterations and \texttt{CPU} solution time for experiments with different nested Krylov methods to solve problem \texttt{QCDx}.}
\vspace{1mm}
\begin{tabular}{rcrcrcrcrcr}
\hline &\multicolumn{2}{c}{\texttt{IT\_in = 150}} &\multicolumn{2}{c}{\texttt{IT\_in = 140}}
&\multicolumn{2}{c}{\texttt{IT\_in = 130}}&\multicolumn{2}{c}{\texttt{IT\_in = 120}}\\
[-2pt]\cmidrule(l{0.7em}r{0.7em}){2-3} \cmidrule(l{0.7em}r{0.6em}){4-5}\cmidrule(l{0.7em}r{0.7em}){6-7}
\cmidrule(l{0.7em}r{0.7em}){8-9} \\[-11pt]
\texttt{Solver} &\texttt{IT\_out} &\texttt{CPU} & $\texttt{IT\_out}$ &\texttt{CPU} &\texttt{IT\_out} &
$\texttt{CPU}$&\texttt{IT\_out} &\texttt{CPU}   \\
\hline
\texttt{Hessen-FCMRH}    &7   &22.73  &7   &20.67   &8   &21.39  &8   &18.95  \\
\texttt{Hessen-FGMRES}   &6   &19.43  &6   &17.52   &7   &18.47  &8   &19.04  \\
\texttt{FOM-FCMRH}       &5   &25.64  &5   &23.04   &6   &24.96  &7   &25.62  \\
\texttt{FOM-FGMRES}      &5   &25.97  &6   &28.46   &6   &25.61  &7   &25.73  \\
\hline
\end{tabular}
\label{tab9xx}
\end{center}
\end{table}

Here problem \texttt{QCDx} is denoted as problem $\Pi_2$ in Section~\ref{sec:5.2}. Experiments with nested iterations using a different number of inner iterations \texttt{IT\_in} are listed in Table \ref{tab9xx}. Our \texttt{Hessen-FGMRES} method with \texttt{IT\_in = 140} has the best overall performance among the four solvers in terms of elapsed CPU time. Again, nested solvers based on the Arnoldi procedure tend to require more CPU time than those based on the Hessenberg procedure. The larger the number of inner iterations, the smaller the number of outer steps, while elapsed CPU time costs do not always follow the same trend. Also in this case we conclude that the \texttt{Hessen-FGMRES} method can be considered as a robust choice for this test
problem. The \texttt{Hessen-FCMRH} method is also an interesting alternative.
%
\section{Conclusions}
\label{sec6}
The paper presents two contributions to the development of shifted Krylov subspace methods built upon
the Hessenberg process for the efficient solution of shifted linear systems. Firstly, we explore the
algorithmic relation between the GMRES($m$) and CMRH($m$) methods. The reduced memory and algorithmic
complexity of CMRH($m$) motivated us to generalize this algorithm for solving shifted linear systems.
The experiments reported in this paper show the effectiveness of the shifted CMRH($m$) method against
\texttt{sGMRES($m$)}, \texttt{sAd-SGMRES($m$)} and \texttt{sQMRIDR($s$)}($s = 1,2$) in terms of elapsed
CPU time. Then, a new nested iterative framework of shifted linear systems is proposed based on the
Hessenberg procedure. More precisely, we proposed three nested iterative solvers: \texttt{Hessen-FCMRH},
\texttt{Hessen-FGMRES} and \texttt{FOM-FCMRH} for shifted linear systems. The first two often converge
significantly faster than \texttt{FOM-FGMRES} introduced in \cite{MBMBG}. In particular, numerical experiments
involving time integration of 3D fractional/partial differential equations are reported to illustrate
the advantages of the proposed nested iterative solvers. We showed that these algorithms can be very
effective to use in numerical schemes that require to evaluate the action of the matrix function on
a vector (see e.g. \cite{JACWL,RGAFMP}) for three-dimensional time-dependent (fractional) convection-diffusion-reaction equations. This point can be regarded as the second
contribution of our manuscript.

In our experiments, the number of inner iterations was often large,
potentially leading to high computational and memory requirements on realistic applications. As an outlook for the future, we plan to test shifted Krylov subspace solvers based on short-term vector recurrences, such as shifted BiCGSTAB($\ell$), shifted BiCRSTAB, shifted IDR($s$), shifted GPBiCG, shifted TFQMR and shifted QMRIDR($s$) as inner solvers in the nested iterative framework for shifted linear systems, see e.g. \cite{MBMBG}. Meanwhile, the restarting technique can often remedy memory problems related to long-term recurrence Krylov subspace methods, e.g., FOM, GMRES, and CMRH. Thus, we are also interested to develop restarted shifted versions of FOM, GMRES, Hessenberg, and CMRH as inner solvers in our future research.
\section*{Acknowledgement}
{\em The authors are greatly grateful to Prof. Gang Wu and Prof. Yan-Fei Jing for their constructive
discussions and insightful comments. Especially, we would like to thank Dr. Mohammed Heyouni and
Dr. Jens-Peter M. Zemke for making their MATLAB package of Ref. \cite{SDMHPM} and codes of the
multi-shifted QMRIDR($s$) method, respectively, available for us. This research is supported by
NSFC (61772003, 11601323, and 11801463), the Fundamental Research Funds for the Central Universities
(JBK1902028) and the Ministry of Education of Humanities and Social Science Layout Project (19JYA790094).}


{\small\begin{thebibliography}{99}
\bibitem{VDLKMZ}
V. Druskin, L. Knizhnerman, M. Zaslavsky, Solution of large scale evolutionary problems
using rational krylov subspaces with optimized shifts, {\em SIAM J. Sci. Comput.}, 31(5)
(2009), pp. 3760-3780.
\bibitem{LNTJAC}
L.N. Trefethen, J.A.C. Weideman, The exponentially convergent trapezoidal rule, {\em SIAM
Rev.}, 56(3) (2014), pp. 385-458.
\bibitem{JACWL}
J.A.C. Weideman, L.N. Trefethen, Parabolic and hyperbolic contours for computing the
Bromwich integral, {\em Math. Comp.}, 76(259) (2007), pp. 1341-1356.
\bibitem{RGAFMP}
R. Garrappa, M. Popolizio, On the use of matrix functions for fractional partial differential
equations, {\em Math. Comput. Simulation}, 81(5) (2011), pp. 1045-1056.
\bibitem{BNDYS}
B.N. Datta, Y. Saad, Arnoldi methods for large Sylvester-like observer matrix equations,
and an associated algorithm for partial spectrum assignment, {\em Linear Algebra Appl.},
154-156 (1991), pp. 225-244.
\bibitem{MIADBS}
M.I. Ahmad, D.B. Szyld, M.B. van Gijzen, Preconditioned multishift BiCG for $\mathcal{H}_2
$-optimal model reduction, {\em SIAM J. Matrix Anal. Appl.}, 38(2) (2017), pp. 401-424.
\bibitem{TITSUN}
T. Sakurai, H. Sugiura, A projection method for generalized eigenvalue problems using
numerical integration, {\em J. Comput. Appl. Math.}, 159(1) (2003), pp. 119-128.
\bibitem{JCRTB}
J.C.R. Bloch, T. Breu, A. Frommer, S. Heybrock, K. Sch\"{a}fer, T. Wettig, Short-recurrence
Krylov subspace methods for the overlap Dirac operator at nonzero chemical potential, {\em
Comput. Phys. Commun.}, 181(8) (2010), pp. 1378-1387.
\bibitem{AKSTB}
A.K. Saibaba, T. Bakhos, P.K. Kitanidis, A flexible Krylov solver for shifted systems with
application to oscillatory hydraulic tomography, {\em SIAM J. Sci. Comput.}, 35(6) (2013),
pp. A3001-A3023.
\bibitem{GWYCWXQJ}
G. Wu, Y.-C. Wang, X.-Q. Jin, A preconditioned and shifted GMRES algorithm for the PageRank
problem with multiple damping factors, {\em SIAM J. Sci. Comput.}, 34(5) (2012), pp. A2558-A2575.
\bibitem{Baumann18}
M. Baumann, M.B. van Gijzen, Convergence and complexity study of GMRES variants for solving
multi-frequency elastic wave propagation problems, \emph{J. Comput. Sci.}, 26 (2018), pp.
285-293.
\bibitem{BJKSS}
B. Jegerlehner, Krylov Space Solvers for Shifted Linear Systems, {\em arXiv preprint},
\href{http://arxiv.org/abs/hep-lat/9612014}{hep-lat/9612014}, December 15, 1996, 16 pages.
\bibitem{AFBCG}
A. Frommer, BiCGStab($\ell$) for families of shifted linear systems, {\em Computing},
70(2) (2003), pp. 87-109.
\bibitem{RWFS}
R.W. Freund, Solution of shifted linear systems by quasi-minimal residual iterations, in
{\em Numerical Linear Algebra: Proceedings of the Conference in Numerical Linear Algebra
and Scientific Computation (L. Reichel, A. Ruttan, R. S. Varga eds.)}, Kent, OH, USA,
March 13-14, 1992, Walter de Gruyter, Berlin, (1993), pp. 101-121.
\bibitem{LDTSSLZ}
L. Du, T. Sogabe, S.-L. Zhang, IDR($s$) for solving shifted nonsymmetric linear systems,
{\em J. Comput. Appl. Math.}, 274 (2015), pp. 35-43.
\bibitem{MBMBG}
M. Baumann, M.B. van Gijzen, Nested Krylov methods for shifted linear systems, {\em SIAM
J. Sci. Comput.}, 37(5) (2015), pp. S90-S112.
\bibitem{MBGGLG}
M.B. van Gijzen, G.L.G. Sleijpen, J.-P. M. Zemke, Flexible and multi-shift induced
dimension reduction algorithms for solving large sparse linear systems, {\em Numer.
Linear Algebra Appl.}, 22(1) (2015), pp. 1-25.
\bibitem{VSRF}
V. Simoncini, Restarted full orthogonalization method for shifted linear systems,
{\em BIT}, 43(2) (2003), pp. 459-466.
\bibitem{YFJTZH}
Y.-F. Jing, T.-Z. Huang, Restarted weighted full orthogonalization method for shifted
linear systems, {\em Comput. Math. Appl.}, 57(9) (2009), pp. 1583-1591.
\bibitem{JFYGJY}
J.-F. Yin, G.-J. Yin, Restarted full orthogonalization method with deflation for shifted
linear systems, {\em Numer. Math. Theor. Meth. Appl.}, 7(3) (2014), pp. 399-412.
\bibitem{XMGTZHG}
X.-M. Gu, T.-Z. Huang, G. Yin, B. Carpentieri, C. Wen, L. Du, Restarted Hessenberg method
for solving shifted nonsymmetric linear systems, {\em J. Comput. Appl. Math.}, 331 (2018),
pp. 166-177.
\bibitem{YSaadI}
Y. Saad, {\em Iterative Methods for Sparse Linear Systems (second ed.)}, SIAM, Philadelphia,
PA, 2003.
\bibitem{TSEX}
T. Sogabe, Extensions of the Conjugate Residual Method (Ph.D dissertation), Department
of Applied Physics, University of Tokyo, Japan, 2006. Available online at
\url{http://www.ist.aichi-pu.ac.jp/person/sogabe/thesis.pdf}.
\bibitem{XMGTZHJ}
X.-M. Gu, T.-Z. Huang, J. Meng, T. Sogabe, H.-B. Li, L. Li, BiCR-type methods for
families of shifted linear systems, {\em Comput. Math. Appl.}, 68(7) (2014), pp. 746-758.
\bibitem{MDRMA}
M. Dehghan, R. Mohammadi-Arani, Generalized product-type methods based on bi-conjugate
gradient (GPBiCG) for solving shifted linear systems, {\em Comput. Appl. Math.}, 36(4)
(2017), pp. 1591-1606.
\bibitem{KAESS}
K. Ahuja, E. de Sturler, S. Gugercin, E.R. Chang, Recycling BiCG with an application
to model reduction, {\em SIAM J. Sci. Comput.}, 34(4) (2012), pp. A1925-A1949.
\bibitem{KAPBES}
K. Ahuja, P. Benner, E. de Sturler, L. Feng, Recycling BiCGSTAB with an application
to parametric model order reduction, {\em SIAM J. Sci. Comput.}, 37(5) (2015), pp.
S429-S446.
\bibitem{JHWK}
J.H. Wilkinson, {\em The Algebraic Eigenvalue Problem}, Clarendon Press, Oxford, UK, 1965.
\bibitem{HSCMRH}
H. Sadok, CMRH: A new method for solving nonsymmetric linear systems based on the
Hessenberg reduction algorithm, {\em Numer. Algorithms}, 20(4) (1999), pp. 303-321.
\bibitem{MHHSA}
M. Heyouni, H. Sadok, A new implementation of the CMRH method for solving dense linear
systems, {\em J. Comput. Appl. Math.}, 213(2) (2008), pp. 387-399.
\bibitem{DS1999}
D. Stephens, ELMRES: An Oblique Projection Method to Solve Sparse Non-Symmetric Linear Systems
(Ph.D Dissertation), Florida Institute of Technology, Melbourne, FL, 1999. Available online
at \url{http://www.ncsu.edu/hpc/Documents/Publications/gary_howell/stephens.ps}.
\bibitem{GHDS1}
G. Howell, D. Stephens, ELMRES, an Oblique Projection Method for Solving Systems of Sparse
Linear Equations, {\em Technical Report}, Florida Institute of Technology, Melbourne, FL,
2000, 18 pages. Available online at \url{http://ncsu.edu/hpc/Documents/Publications/gary_howell/elmres320.ps}.
\bibitem{MHMDHG}
M. Heyouni, M\'{e}thode de Hessenberg G\'{e}n\'{e}ralis\'{e}e et Applications (Ph.D thesis),
Universit\'{e} des Sciences et Technologies de Lille, D\'{e}cembre, France, 1996. Available
online at \url{http://ori.univ-lille1.fr/notice/view/univ-lille1-ori-134864}. (in Franch)
\bibitem{KZCG}
K. Zhang, C. Gu, A flexible CMRH algorithm for nonsymmetric linear systems, {\em J. Appl.
Math. Comput.}, 45(1-2) (2014), pp. 43-61.
\bibitem{SDMHPM}
S. Duminil, M. Heyouni, P. Marion, H. Sadok, Algorithms for the CMRH method for dense linear
systems, {\em Numer. Algorithms}, 71(2) (2016), pp. 383-394.
\bibitem{Lai2011}
J. Lai, L. Lu, S. Xu, A polynomial preconditioner for the CMRH algorithm, {\em Math.
Probl. Eng.}, 2011 (2011), Article ID 545470, 12 pages. DOI:
\href{http://dx.doi.org/10.1155/2011/545470}{10.1155/2011/545470}.
\bibitem{SDAP}
S. Duminil, A parallel implementation of the CMRH method for dense linear systems,
{\em Numer. Algorithms}, 63(1) (2013), pp. 127-142.
\bibitem{YSMHG}
Y. Saad, M.H. Schultz, GMRES: A generalized minimal residual algorithm for solving
nonsymmetric linear systems, {\em SIAM J. Sci. Stat. Comput.}, 7(3) (1986), pp. 856-869.
\bibitem{AFUG}
A. Frommer, U. Gl\"{a}ssner, Restarted GMRES for shifted linear systems, {\em SIAM J.
Sci. Comput.}, 19(1) (1998), pp. 15-26.
\bibitem{GGRG}
G. Gu, Restarted GMRES augmented with harmonic Ritz vectors for shifted linear systems,
{\em Int. J. Comput. Math.}, 82(7) (2005), pp. 837-849.
\bibitem{DDRBMW}
D. Darnell, R.B. Morgan, W. Wilcox, Deflated GMRES for systems with multiple shifts
and multiple right-hand sides, {\em Linear Algebra Appl.}, 429(10) (2008), pp. 2415-2434.
\bibitem{KMSDBS}
K.M. Soodhalter, D.B. Szyld, F. Xue, Krylov subspace recycling for sequences of shifted
linear systems, {\em Appl. Numer. Math.}, 81 (2014), pp. 105-118.
\bibitem{KMST}
K.M. Soodhalter, Two recursive GMRES-type methods for shifted linear systems with general
preconditioning, {\em Electron. Trans. Numer. Anal.}, 45 (2016), pp. 499-523.
\bibitem{Jing2017}
Y.-F. Jing, P. Yuan, T.-Z. Huang, A simpler GMRES and its adaptive variant for shifted
linear systems, {\em Numer. Linear Algebra Appl.}, 24(1) (2017), Article No. e2076, 7
pages. DOI: \href{http://dx.doi.org/10.1002/nla.2076}{10.1002/nla.2076}.
\bibitem{HSDBS}
H. Sadok, D.B. Szyld, A new look at CMRH and its relation to GMRES, {\em BIT}, 52(2)
(2012), pp. 485-501.
\bibitem{Meurant17}
G. Meurant, An optimal Q-OR Krylov subspace method for solving linear systems, {\em
Electron. Trans. Numer. Anal.}, 47 (2017), pp. 127-152.
\bibitem{VSDBS}
V. Simoncini, D.B. Szyld, The effect of non-optimal bases on the convergence of Krylov
subspace methods, {\em Numer. Math.}, 100(4) (2005), pp. 711-733.
\bibitem{JDTGM}
J. Duintjer Tebbens, G. Meurant, On the convergence of Q-OR and Q-MR Krylov methods for
solving nonsymmetric linear systems, {\em BIT}, 56(1) (2016), pp. 77-97.
\bibitem{Imakura13}
A. Imakura, T. Sogabe, S.-L. Zhang, An efficient variant of the restarted shifted GMRES
for solving shifted linear systems, {\em J. Math. Res. Appl.}, 33(2) (2013), pp. 127-141.
\bibitem{TSTHSLZ}
T. Sogabe, T. Hoshi, S.-L. Zhang, T. Fujiwara, A numerical method for calculating the Green's
function arising from electronic structure theory, in \textit{Frontiers of Computational
Science}, Y. Kaneda, H. Kawamura, M. Sasai (Eds.), Springer-Verlag, Berlin/Heidelberg, (2007),
pp. 189-195. DOI: \href{https://doi.org/10.1007/978-3-540-46375-7_24}{10.1007/978-3-540-46375-7\_24}.
\bibitem{GGXZLL}
G.-D. Gu, X.-L. Zhou, L. Lin, A flexible preconditioned Arnoldi method for shifted linear
systems, {\em J. Comp. Math.}, 25(5) (2007), pp. 522-530.
\bibitem{YSaadA}
Y. Saad. A flexible inner-outer preconditioned GMRES algorithm, {\em SIAM J. Sci. Comput.},
14(2) (1993), pp. 461-469.
\bibitem{VSDBSF}
V. Simoncini, D.B. Szyld, Flexible inner-outer Krylov subspace methods, {\em SIAM J. Numer.
Anal.}, 40(6) (2003), pp. 2219-2239.
\bibitem{BCIDLGGS}
B. Carpentieri, I.S. Duff, L. Giraud, G. Sylvand, Combining fast multipole techniques and an
approximate inverse preconditioner for large electromagnetism calculations, {\em SIAM J. Sci.
Comput.}, 27(3) (2005), pp. 774--792.
\bibitem{DSTHBCYJ}
D.-L. Sun, T.-Z. Huang, B. Carpentieri, Y.-F. Jing, Flexible and deflated
variants of the block shifted {GMRES} method, {\em J. Comput. Appl. Math.},
345 (2019), pp. 168-183.
\bibitem{TDYHT}
T. Davis, Y. Hu, The University of Florida Sparse Matrix Collection, {\em ACM Trans. Math.
Software}, 38(1) (2011), Article No. 1, 25 pages. Avaiable online at \url{https://sparse.tamu.edu/}.
\bibitem{IPFD}
I. Podlubny, {\em Fractional Differential Equations}, Academic Press, San Diego, CA, 1999.
\bibitem{MCCES}
S. Zhai, D. Gui, P. Huang, X. Feng, A novel high-order ADI method for 3D fractional
convection-diffusion equations, {\em Int. Commun. Heat Mass Transfer}, 66 (2015), pp.
212-217.
\bibitem{LNTJAC2}
L.N. Trefethen, J.A.C. Weideman, T. Schmelzer, Talbot quadratures and rational approximations,
{\em BIT}, 46(3) (2006), pp. 653-670.
\bibitem{MCPKAO}
M. Caliari, P. Kandolf, A. Ostermann, S. Rainer, Comparison of software for computing the
action of the matrix exponential, {\em BIT}, 54(1) (2014), pp. 113-128.
\end{thebibliography}}
\end{document}